%% file: arxiv.tex
\pdfoutput=1
\documentclass{amsart}
\usepackage[left=3.25cm, right=3.25cm, top=3cm, bottom=3cm]{geometry}

\usepackage{lipsum}
\usepackage{amsfonts}
\usepackage{graphicx}
\usepackage{float}
\usepackage{algorithm}
\usepackage{algorithmic}
\usepackage{multirow,makecell}
\usepackage{bm,amsfonts,booktabs,amssymb}
\usepackage{amsmath}
\makeatletter
\g@addto@macro\normalsize{%
  \setlength\abovedisplayskip{7pt}
  \setlength\belowdisplayskip{7pt}
  \setlength\abovedisplayshortskip{8pt}
  \setlength\belowdisplayshortskip{8pt}
}
\makeatother
\usepackage{mathtools}
\usepackage{enumitem}
\usepackage{subcaption}
\usepackage[font={small,it}]{caption}
\usepackage{comment}
\usepackage{nicematrix}
\usepackage[hidelinks]{hyperref}
\usepackage{cleveref}
\usepackage{autonum}
\usepackage{longtable}
\usepackage{amsopn}
\usepackage{tikz}
\usetikzlibrary{positioning}
\usetikzlibrary{arrows}
\usepackage{fancyvrb} 
\VerbatimFootnotes    

\ifpdf
  \DeclareGraphicsExtensions{.eps,.pdf,.png,.jpg}
\else
  \DeclareGraphicsExtensions{.eps}
\fi

\newcommand{\refassump}[2]{\hyperref[#1]{#2}}

\newcolumntype{P}[1]{>{\centering\arraybackslash}p{#1}}
\newcolumntype{M}[1]{>{\centering\arraybackslash}m{#1}}

\newtheorem{theorem}{Theorem}[section]
\newtheorem{corollary}{Corollary}[section]
\newtheorem{lemma}{Lemma}[section]

\newtheorem{remark}{Remark}[section]

\everypar{\looseness=-1}

\newcommand{\kkeywords}[1][]{
\smallskip
\noindent \textbf{Keywords.} #1
}

\crefname{equation}{}{}

\newcommand{\textalgo}[1]{\texttt{\textbf{#1}}}

\begin{document}
\renewcommand{\thefootnote}{\fnsymbol{footnote}}
\footnotetext[1]{equal contribution}
\footnotetext[2]{Program in Applied and Computational Mathematics, Princeton University (dhruv.kohli@princeton.edu)}
\footnotetext[3]{Halicio\u{g}lu Data Science Institute, UC San Diego (jeh020@ucsd.edu, chholtz@ucsd.edu, gmishne@ucsd.edu)}
\footnotetext[4]{Department of Mathematics, UC San Diego (acloninger@ucsd.edu)}

\title[Robust boundary detection and density estimation using doubly stochastic kernel]{Robust boundary detection and density estimation using doubly stochastic scaling of the Gaussian kernel}
\author[D. Kohli, J. He, C. Holtz, A. Cloninger, G. Mishne]{Dhruv Kohli${}^{\dagger\ast}$, Jesse He${}^{\ddagger\ast}$, Chester Holtz${}^\ddagger$, Alexander Cloninger${}^{\S\ddagger}$,  Gal Mishne${}^\ddagger$}

\renewcommand{\thefootnote}{\arabic{footnote}}

\begin{abstract}
This paper addresses the problem of detecting boundary points and estimating the sampling density of a dataset derived from a compact manifold with boundary, potentially in the presence of noise.
We extend recent advances in doubly stochastic scaling of the Gaussian heat kernel via Sinkhorn iterations to this setting.
Our main contributions are: (a) deriving a characterization of the scaling factors for manifolds with boundary, (b) developing a boundary direction estimator aimed at identifying boundary points followed by a boundary-corrected kernel density estimator based on doubly stochastic kernel and local principal component analysis, and (c) demonstrating through simulations that the resulting estimates of the boundary points and the sampling density outperform the standard Gaussian kernel-based approach, particularly under noisy conditions.

\kkeywords{Manifold with boundary, doubly stochastic kernel, boundary detection, kernel density estimation.}

\end{abstract}

\maketitle

\section{Introduction}
\input{sections/intro}

\section{Characterizing doubly stochastic scaling factors on manifolds with boundary}
\label{sec:main}
\input{sections/theory}

\section{Estimation of boundary and sampling density using doubly stochastic kernel}
\label{sec:algo}
\input{sections/algo}

\subsection{Density estimation using doubly stochastic kernel}
\label{subsec:ds_kde}
\input{sections/kde}

\section{Experiments}
\label{sec:experiments}
\input{sections/exp}

\section{Conclusion}
\label{sec:conclusion}
\input{sections/conc}

\appendix
\section{Proofs}
\label{sec:proofs}
\input{sections/proofs}

\section*{Acknowledgments}
This work was partially funded by NSF award 2217058. CH also wishes to acknowledge an award from the W.M. Keck foundation. 

\bibliographystyle{siamplain}
\bibliography{main}

\end{document}

%% file: sections/intro.tex
We address the problem of identifying points on or near the boundary of a compact manifold~\cite{lee2018introduction} using a potentially noisy set of sampled points.
The knowledge of such points is crucial for numerically solving partial differential equations with boundary conditions~\cite{VAUGHN2024101593}, unbiased kernel density estimation~\cite{schuster1985incorporating,jones1996simple,karunamuni2005boundary}, constructing low-distortion local eigenmaps for regions near the boundary~\cite{jones2007universal,ldle1v2,ratsv2}, and evaluating the performance of clustering algorithms~\cite{jung2018network} in preserving the structure within clusters.

Estimating the boundary of a data manifold remains a relatively underexplored problem. Berry and Sauer~\cite{berry2017density} proposed an estimation strategy aimed at the general case of data sampled from compact Riemannian manifolds in arbitrary dimensions. They employed the standard Gaussian kernel density estimator (KDE) along with a boundary direction estimator (BDE) to numerically solve for the distances of the points from the boundary. Similarly, Calder, Park, and Slep\v{c}ev~\cite{calder2022boundary} introduced a boundary estimator that replaces the Gaussian kernel in the KDE and BDE with a binary cutoff kernel. In this work, we show that boundary estimates based on both the standard Gaussian kernel and the binary cutoff kernel are highly sensitive to noise. Moreover, the latter method is also sensitive to the curvature of the underlying manifold, as demonstrated in our experiments.

A recent line of research~\cite{marshall2019manifold,landa2021doubly,landa2023robust,cheng2024bi} investigated the concept of doubly stochastic scaling via Sinkhorn iterations~\cite{cuturi2013sinkhorn} applied to the standard Gaussian heat kernel.
These studies demonstrated that the doubly stochastic kernel leads to significant improvements in kernel density estimation on a closed manifold, and in approximating the eigenvectors of the manifold Laplacian, especially in the presence of high-dimensional heteroskedastic and outlier-type noise. Additionally, Cheng and Landa~\cite{cheng2024bi} established convergence results for the bi-stochastically normalized Laplacian to a weighted manifold Laplacian and also demonstrated its resilience to outlier noise. In this paper, we follow a similar approach by utilizing the doubly stochastic kernel to obtain robust boundary and density estimates for data points sampled from a manifold with noise. The contributions of our paper are as follows.

First, building on the characterization by Landa and Cheng~\cite{landa2023robust} of the scaling factors $\rho_h(\mathbf{x})$ that transform a Gaussian heat kernel $k_{h}(\mathbf{x},\mathbf{y})$ with bandwidth $h$ into a doubly stochastic kernel on a \emph{closed manifold}, we derive an analogous characterization of the scaling factors in the setting of a \emph{manifold with boundary}. Our second-order characterization, presented in~\Cref{thm:rho_hx_char}, takes the form:
\begin{equation}
    A_h(\mathbf{x})\rho_h(\mathbf{x})^4q(\mathbf{x})^2 - B_h(\mathbf{x}))\rho_h(\mathbf{x})^2q(\mathbf{x}) + \pi^{d} = \mathcal{O}(h^2)
\end{equation}
where $d$ is the dimension of the manifold. The coefficients $A_h(\mathbf{x})$ and $B_h(\mathbf{x})$ depend on the distance $b_{\mathbf{x}}$ of $\mathbf{x}$ from the boundary, the mean curvature induced by the boundary at $\mathbf{x}$, and the derivative of the sampling density $q(\mathbf{x})$ in the direction normal to the boundary at $\mathbf{x}$. Since these quantities are typically unknown for a given data manifold, for practical purposes, we obtain a relaxed first-order characterization in~\Cref{corollary:rho_h_first_order_char},
\begin{equation}
    \rho_h(\mathbf{x})^2q(\mathbf{x}) = \zeta_h(\mathbf{x}) +\mathcal{O}(h) \label{eq:rho_h_first_order_char}
\end{equation}
where $\zeta_h(\mathbf{x})$ depends only on the distance of $\mathbf{x}$ from the boundary, i.e. on $b_{\mathbf{x}}$.

Second, given a noisy point cloud $\{\widetilde{\mathbf{x}}_i\}_1^n$, we introduce a BDE $\{\boldsymbol{\nu}_i\}_1^n$ based on the doubly stochastic heat kernel $\mathbf{W}$~\cite{marshall2019manifold,landa2021doubly,landa2023robust,cheng2024bi} and local principal component analysis (PCA)~\cite{cleveland1988locally,arias2011spectral,singer2012vector,tyagi2013tangent,cheng2013local,kaslovsky2014non,yu2015useful,mohammed2017manifold,aamari2019nonasymptotic,gilbert2025pca}. Specifically,
\begin{equation}
    \boldsymbol{\nu}_i \coloneqq \frac{1}{n-1}\sum_{j=1}^{n}\mathbf{W}_{ij} \widetilde{\mathbf{U}}_{i}^T(\widetilde{\mathbf{x}}_j-\widetilde{\mathbf{x}}_i) \label{eq:nu_teaser}
\end{equation}
where the columns of $\widetilde{\mathbf{U}}_{i} \in \mathbb{R}^{m \times d}$ represent the $d$ principal directions obtained by applying PCA in a neighborhood of the $i$-th noisy data point $\widetilde{\mathbf{x}}_i$.
In~\Cref{theorem:Enu}, we establish the convergence of $\boldsymbol{\nu}_i$ to the normal direction to the boundary $\boldsymbol{\eta}_{\mathbf{x}_i}$ at $\mathbf{x}_i$, the clean counterpart of $\widetilde{\mathbf{x}}_i$, up to a global scaling factor that depends only on the distance to the boundary. To keep the focus of the paper on doubly stochastic scaling, our convergence analysis assumes access to the projector onto the tangent space. In our experiments, this projector is approximated via local PCA on noisy data points. Under this assumption, we show that the following holds with high probability:
\begin{equation}
    \boldsymbol{\nu}_i \xrightarrow[]{\ n \rightarrow \infty \ } h\beta_h(\mathbf{x}_i) \boldsymbol{\eta}_{\mathbf{x}_i}  + \mathcal{O}\left(h^2\right).\label{eq:E[nu]_teaser}
\end{equation}
where $\beta_h(\mathbf{x})$ depends solely on the distance $b_{\mathbf{x}}$ of $\mathbf{x}$ from the boundary. Moreover, $|\beta_h(\mathbf{x})|$ is a strictly decreasing function of $b_{\mathbf{x}}$ which allows boundary points to be identified by simple thresholding of the normed\footnote{$\left\|\cdot \right\|$ corresponds to the $2$-norm throughout this work.} vectors $\left\|\boldsymbol{\nu}_i\right\|$. 
Our experimental results show that this method results in a significant improvement in boundary estimation compared to the standard Gaussian kernel-based approach~\cite{berry2017density}, particularly in the presence of heteroskedastic noise.
    
Finally, we propose a two-parameter family of kernel density estimators for manifolds with boundary that corrects the boundary-induced bias inherent in standard KDEs~\cite{marron1994transformations,cowling1996pseudodata,karunamuni2005boundary} by explicitly incorporating distances to the boundary. This approach extends the robust KDE framework of~\cite{landa2023robust} from closed manifolds to manifolds with boundary while remaining effective under heteroskedastic noise.

The organization of the paper is as follows: in~\Cref{subsec:setup} we introduce the setup. In~\Cref{sec:main}, we provide a characterization of the doubly stochastic scaling factors on manifolds with boundary. In~\Cref{sec:algo}, we present our BDE and KDE, and analyze their convergence. In~\Cref{sec:experiments}, we showcase our experimental results demonstrating the effectiveness of our method, followed by conclusions in~\Cref{sec:conclusion}. Proofs of our results are provided in~\Cref{sec:proofs}.

\section{Setup and Preliminaries}
\label{subsec:setup}
Let $\mathcal{M}$ be a $d$-dimensional compact Riemannian manifold with boundary $\partial \mathcal{M}$ embedded in an $m$-dimensional Euclidean space.
Let $|\mathcal{M}|$ be the volume of $\mathcal{M}$, $\boldsymbol{\eta}_{\mathbf{x}}$ be the normal direction at $\mathbf{x}$ to the boundary, $H(\mathbf{x})$ be the mean curvature of the hypersurface parallel to $\partial \mathcal{M}$ intersecting $\mathbf{x}$ (which depends on the second fundamental form of $\partial \mathcal{M} \subset \mathcal{M}$), $b_{\mathbf{x}}$ denote the distance of $\mathbf{x}$ from the boundary and $dV(\mathbf{x})$ be the volume form of $\mathcal{M}$ at $\mathbf{x}$.

Suppose $(\mathbf{x}_i)_1^n \subset \mathbb{R}^m$ are sampled independently and identically from a smooth and positive probability distribution $q(\mathbf{x})$ on $\mathcal{M}$.
Let $\widetilde{\mathbf{x}}_i = \mathbf{x}_i + \boldsymbol{\varepsilon}(\mathbf{x}_i)$ be the observed noisy samples where $\{\boldsymbol{\varepsilon}(\mathbf{x}_i)\}_1^n \subset \mathbb{R}^m$ are independent, perhaps non-identical, sub-Gaussian random variables with zero mean and a sub-Gaussian norm of $\left\|\boldsymbol{\varepsilon}(\mathbf{x}_i)\right\|_{\psi_2}$. The noiseless setting is characterized by the constraint $\varepsilon(\mathbf{x}_i) = 0$ for all $i \in {1,\ldots,n}$.
For $h > 0$, define the Gaussian kernel $k_h$ and the affinity matrix $\mathbf{K}_{ij}$ by
\begin{align}
    k_h(\widetilde{\mathbf{x}}_i, \widetilde{\mathbf{x}}_j) &\coloneqq e^{-\left\|\widetilde{\mathbf{x}}_i-\widetilde{\mathbf{x}}_j\right\|^2/h^2},\\
    \mathbf{K}_{ij} &\coloneqq \begin{cases}k_h(\widetilde{\mathbf{x}}_i, \widetilde{\mathbf{x}}_j), & i \neq j\\ 0, & i = j. \label{eq:Kij}
    \end{cases}
\end{align}
Given the kernel matrix $\mathbf{K}$, the doubly stochastic matrix $\mathbf{W} \in \mathbb{R}^{n \times n}$ is defined as in~\cite{marshall2019manifold,landa2021doubly} by
\begin{equation}
    \mathbf{W}_{ij} = \mathbf{d}_i\mathbf{K}_{ij}\mathbf{d}_j \text{ and }\sum_{j=1}^{n}\mathbf{W}_{ij} = 1, \label{eq:W}
\end{equation}
where $\mathbf{d} \in \mathbb{R}^n$ are scaling factors. The existence of $\mathbf{d}$ follows from \cite{landa2021doubly} and these are obtained efficiently using Sinkhorn iterations~\cite{cuturi2013sinkhorn}. In the continuous setting, these correspond to a positive function $\rho_{h}$ on $\mathcal{M}$ such that
\begin{equation}
    \frac{1}{\pi^{d/2}h^d}\int_{\mathcal{M}}\rho_{h}(\mathbf{x})k_{h}(\mathbf{x},\mathbf{y})\rho_{h}(\mathbf{y}) q(\mathbf{y})dV(\mathbf{y}) = 1. \label{eq:rho_constraint}
\end{equation}
The existence of a smooth $\rho_h$ on $\mathcal{M}$ that satisfies the above constraint follows from~\cite{marshall2019manifold,knopp1968note,borwein1994entropy}. In this work, we only need the following assumption:
\begin{enumerate}
    \item[\textit{(A1)}]\label{assump:rho} There exists a positive scaling function $\rho_h \in \mathcal{C}^{3}(\mathcal{M})$ that satisfies Eq.~\cref{eq:rho_constraint}. 
\end{enumerate}

Under the sub-Gaussian noise setup, Landa and Cheng in~\cite{landa2023robust} established a connection between the doubly stochastic kernel obtained from noisy observations and the Gaussian kernel obtained from clean samples. We use their result and the underlying assumptions, as stated below, in the convergence analysis of our BDE in~\Cref{sec:algo}.
\begin{lemma}[Theorem 2.1 in \cite{landa2023robust}]
\label{lem:setup}
Under the assumptions:
\begin{enumerate}
    \item[(A2)]\label{assump:q_old} The sampling density $q$ is positive and continuous on $\mathcal{M}$. 
    \item[(A3)]\label{assump:manifold} $\left\|\mathbf{x}\right\| \leq 1$ for all $\mathbf{x} \in \mathcal{M}$.
    \item[(A4)]\label{assump:noise} There exists $C > 0$ such that $\max_{\mathbf{x} \in \mathcal{M}}\left\|\boldsymbol{\varepsilon}(\mathbf{x}_i)\right\|_{\psi_2} \leq C/(m^{1/4}\sqrt{\log m})$ for all $\mathbf{x} \in \mathcal{M}$.
    \item[(A5)]\label{assump:dim} There exists $\gamma > 0$ such that $m \geq n^{\gamma}$. 
\end{enumerate}
the following holds,
\begin{align}
    \mathbf{W}_{i,j} &= \frac{\rho_{h}(\mathbf{x}_i)k_{h}(\mathbf{x}_i, \mathbf{x}_j)\rho_{h}(\mathbf{x}_j)}{(n-1)\pi^{d/2}h^d}\left(1+\mathcal{E}_{i,j}\right)\label{eq:Wij_result}\\
    \mathbf{d}_{i} &= \frac{\rho_{h}(\mathbf{x}_i)}{\sqrt{(n-1)\pi^{d/2}h^d}}e^{\frac{\left\|\varepsilon(\mathbf{x}_i)\right\|^2}{h^2}}\left(1+\mathcal{E}_{i,i}\right) \label{eq:di}
\end{align}
where $\mathcal{E}_{i,j} = \mathcal{O}^{(h)}_{m,n}(g(m,n))$. In the noiseless case, $g(m,n) = g_0(n)$, whereas in the presence of noise, $g(m,n) = \max\{g_0(n), g_1(m)\}$. Here,
\begin{align}
    g_0(n) &\coloneqq \sqrt{\log n/n}, \label{eq:E_n}\\
    g_1(m) &\coloneqq \max\{m^{-1/4},1/\sqrt{\log m}\} ,\label{eq:E_m} 
\end{align}
and as defined in \cite{landa2023robust}, a random variable $X = \mathcal{O}^{(h)}_{m,n}(g(m,n))$ if there exist $t_0$, $n_0(h)$, $m_0(h)$, $C(h) > 0$ such that for all $n > n_0(h)$, $m \geq m_0(h)$, $$\mathrm{Pr}\{|X| \leq t C(h)g(m,n)\} \geq 1-n^{-t} \text{ for all } t \geq t_0.$$ Moreover,
\begin{enumerate}[leftmargin=*,label=(\roman*)]
    \item if $Y = f^{(h)}(X)$ where $f^{(h)} \in \mathcal{C}^1(\mathbb{R})$ for every $h > 0$ then 
    \begin{equation}
        Y = f^{(h)}(0) + \mathcal{O}^{(h)}_{m,n}(g(m,n)). \label{eq:g^h(X)}
    \end{equation} 
    \item if $X_i = \mathcal{O}^{(h)}_{m,n}(g(m,n))$ for all $i \in [1,n]$ then 
    \begin{equation}
        \max_1^n|X_i| = \mathcal{O}^{(h)}_{m,n}(g(m,n)). \label{eq:maxXi}
    \end{equation}
\end{enumerate}
\end{lemma}

\begin{remark}
    Assumption \refassump{assump:rho}{(A4)} bounds the noise in terms of the ambient dimension while the assumption \refassump{assump:rho}{(A5)} controls the growth of the ambient dimension with the sample size. Overall, the two assumptions facilitate the analysis of doubly stochastic scaling in the large sample size and high dimension regime.
\end{remark}
\begin{remark}
    The above result holds for continuous density $q$ on $\mathcal{M}$. However, we make the following stronger assumption throughout this work.
    \begin{enumerate}
        \item[(A2')]\label{assump:q} The sampling density $q$ is positive and $q \in \mathcal{C}^{3}(\mathcal{M})$.
    \end{enumerate}
\end{remark}

%% file: sections/theory.tex
In~\cite{landa2023robust}, the authors characterized the scaling factors $\rho_h(\mathbf{x})$ as defined in Eq.~(\ref{eq:rho_constraint}) in the setting of a closed manifold. Here, we obtain a characterization of $\rho_h(\mathbf{x})$ on manifolds with boundary. Specifically, we obtain a second order characterization of $\rho_h(\mathbf{x})$ that depends on the sampling density $q(\mathbf{x})$, the mean curvature $H(\mathbf{x})$ induced by the boundary and the distance $b_{\mathbf{x}}$ to the boundary. Since, for a given dataset, the curvature $H(\mathbf{x})$ is not known \emph{a priori}, we further obtain a first order characterization of $\rho_h(\mathbf{x})$ that only depends on $q(\mathbf{x})$ and $b_{\mathbf{x}}$ and later utilize this characterization to obtain robust estimate of points on the boundary.

First, we need the following lemmas: (a) on the third-order approximation of the integral against Gaussian kernel on a manifold with boundary~\cite{coifman2006diffusion,vaughn2019diffusion} as stated below, and (b) on the directional derivative of the projection operator onto the tangent space (see~\Cref{lemma:Px_eqs} in the Appendix).

\begin{lemma}[Theorem~4.7 in \cite{vaughn2019diffusion}]
\label{lem:integral_khf}
Let $\mathcal{M}$ be a compact $d$-dimensional $\mathcal{C}^3$ Riemannian manifold with a $\mathcal{C}^3$ boundary and a normal collar. Then for all $f \in \mathcal{C}^{3}(\mathcal{M})$, $\mathbf{x} \in \mathcal{M}$ and $h$ sufficiently small,
{\normalsize
\begin{align}
    &\frac{1}{h^{d}}\int_{\mathcal{M}} k_{h}(\mathbf{x},\mathbf{y})f(\mathbf{y}) dV(\mathbf{y})\\
    &\hspace{0.75cm} = m^{(0)}_h(\mathbf{x})f(\mathbf{x}) +  hm^{(1)}_{h}(\mathbf{x})\left(\partial_{\boldsymbol{\eta}_{\mathbf{x}}} f(\mathbf{x}) + \frac{d-1}{2}H(\mathbf{x})f(\mathbf{x})\right) + \frac{h^2}{2}\left\{\frac{\pi^{d/2}}{2}\widetilde{\omega}(\mathbf{x})f(\mathbf{x}) + \right.\\
    &\hspace{1.5cm} \left. \frac{\pi^{d/2}}{2}(\Delta f(\mathbf{x}) - \partial_{\boldsymbol{\eta}_{\mathbf{x}}}^{2}f(\mathbf{x})) + m^{(2)}_h(\mathbf{x})(\partial_{\boldsymbol{\eta}_{\mathbf{x}}}^{2}f(\mathbf{x}) - (d-1)H(\mathbf{x}) \partial_{\eta_\mathbf{x}}f(\mathbf{x}))\right\} + \mathcal{O}(h^3)
\end{align}}%
where $H(\mathbf{x})$ is the mean curvature of the hypersurface parallel to $\partial \mathcal{M}$ intersecting $\mathbf{x}$ (which depends on the second fundamental form of $\partial \mathcal{M} \subset \mathcal{M}$), $\widetilde{\omega}(\mathbf{x})$ depends on the scalar curvature of $\mathcal{M}$ and the boundary, and
\begin{align}
    m_{h}^{(0)}(\mathbf{x}) &\coloneqq \frac{\pi^{d/2}}{2}\left(1+\mathrm{erf}\left(\frac{b_{\mathbf{x}}}{h}\right)\right) \label{eq:m0}\\
    m_{h}^{(1)}(\mathbf{x}) &\coloneqq -\frac{\pi^{(d-1)/2}}{2}e^{-b_{\mathbf{x}}^2/h^2}\label{eq:m1}\\
    m^{(2)}_{h}(\mathbf{x}) &\coloneqq \frac{b_{\mathbf{x}}}{h} m^{(1)}_{h}(\mathbf{x}) + \frac{m^{(0)}_{h}(\mathbf{x})}{2}.\label{eq:m2}
\end{align}
\end{lemma}

In the following result, by applying~\Cref{lem:integral_khf} and~\Cref{lemma:Px_eqs} to Eq.~\cref{eq:rho_constraint}, we derive the directional derivatives of the scaling factors $\rho_h(\mathbf{x})$ and a quadratic equation in $\rho_h^2(\mathbf{x})$ that characterizes the scaling factors on a manifold with boundary. These derivatives are utilized in Section~\ref{sec:algo} to establish the convergence of our BDE to the normal direction to the boundary, up to a scaling factor that depends on the distances to the boundary.
\begin{theorem}
\label{thm:rho_hx_char}
Let $\mathcal{M}$ be a compact $d$-dimensional $\mathcal{C}^3$ Riemannian manifold with a $\mathcal{C}^3$ boundary and a normal collar. Suppose assumptions \refassump{assump:rho}{(A1)} and \refassump{assump:q}{(A2')} hold: $\rho_{h} \in \mathcal{C}^3(\mathcal{M})$ be a positive scaling function satisfying Eq.~(\ref{eq:rho_constraint}) and $q \in \mathcal{C}^3(\mathcal{M})$ be the positive sampling density. Then the following holds.

\noindent (I) The derivative of $\rho_{h}$ in the direction of $\boldsymbol{\eta}_{\mathbf{x}}$ is,
\begin{equation}
    h\partial_{\boldsymbol{\eta}_{\mathbf{x}}}\rho_h(\mathbf{x}) = -\frac{2\rho_h(\mathbf{x})^3q(\mathbf{x})\left\{m^{(1)}_h(\mathbf{x}) + hm^{(2)}_h(\mathbf{x})\kappa_{-}(\mathbf{x})\right\}}{\pi^{d/2}+2m^{(2)}_h(\mathbf{x})\rho_h(\mathbf{x})^2 q(\mathbf{x})} + \mathcal{O}(h^2). \label{eq:partial_eta_rho_h2}
\end{equation}
Also, if $\mathbf{v}_{\mathbf{x}} \in T_{\mathbf{x}}\mathcal{M}$ such that $\mathbf{v}_\mathbf{x} \perp \boldsymbol{\eta}_{\mathbf{x}}$ then the derivative of $\rho_{h}$ in the direction of $\mathbf{v}_{\mathbf{x}}$ is,
\begin{equation}
    \partial_{\mathbf{v}_{\mathbf{x}}}\rho_h(\mathbf{x}) = -\frac{\rho_h(\mathbf{x})^3\partial_{\mathbf{v}_{\mathbf{x}}}q(\mathbf{x})}{(1+\rho_h(\mathbf{x})^2 q(\mathbf{x}))} 
     + \mathcal{O}(h). \label{eq:partial_v_rho_h}
\end{equation}

\noindent (II) $\rho_h(\mathbf{x})$ satisfies the equation,
\begin{equation}
    A_h(\mathbf{x})\rho_h(\mathbf{x})^4q(\mathbf{x})^2 - B_h(\mathbf{x})\rho_h(\mathbf{x})^2q(\mathbf{x}) + \pi^{d} = \mathcal{O}(h^2) \label{eq:rho_h_quadratic}
\end{equation}
where $A_h$ and $B_h$ are defined as follows:
\begin{align}
    A_h(\mathbf{x}) &= 2\left\{m^{(1)}_{h}(\mathbf{x})^2 - m^{(2)}_{h}(\mathbf{x})\left(m^{(0)}_{h}(\mathbf{x}) + h(d-1) m^{(1)}_{h}(\mathbf{x})H(\mathbf{x})\right) \right\} \text{ and}\\
    B_h(\mathbf{x}) &= \pi^{d/2}m^{(1)}_{h}(\mathbf{x})\left(h\kappa_{+}(\mathbf{x}) - \frac{2b_{\mathbf{x}}}{h}\right),
\end{align}
where $\kappa_{\pm}(\mathbf{x}) \coloneqq \frac{\partial_{\boldsymbol{\eta}_{\mathbf{x}}}q(\mathbf{x})}{q(\mathbf{x})} \pm \frac{(d-1)}{2}H(\mathbf{x})$ depends on the sampling density and the mean curvature $H(\mathbf{x})$ induced by the boundary at $\mathbf{x}$.
Since $\mathcal{M}$ is compact, the constant in $\mathcal{O}(h^2)$ can be uniformly bounded over the manifold.
\end{theorem}

All the corollaries presented below are extensions of~\cref{thm:rho_hx_char}. By specializing the proof of this theorem to closed and compact manifolds, we obtain the following corollary, which is reminiscent of the characterization of scaling factors presented in \cite{landa2023robust}.
\begin{corollary}
\label{corollary:closed_manifold}
If $\partial\mathcal{M} = \emptyset$ then

$$\rho_h(\mathbf{x})^2q(\mathbf{x}) = 1+\mathcal{O}(h^2)$$
where the constant in $\mathcal{O}(h^2)$ which depends on $\mathbf{x}$ can be uniformly bounded since $\mathcal{M}$ is compact.
\end{corollary}

In fact, for a sufficiently small $h$ the above result holds at the interior points of a manifold with nonempty boundary whenever $b_{\mathbf{x}}/h \gg 0$. The key component of the proof is the asymptotic expansion of the heat kernel in the interior of the manifold~\cite{vaughn2019diffusion, coifman2006diffusion}, combined with the fact that when $b_{\mathbf{x}}/h \rightarrow \infty$ then $m^{(0)}_{h}(\mathbf{x}) \rightarrow \pi^{d/2}$, $\frac{b_{\mathbf{x}}}{h}m^{(1)}_{h}(\mathbf{x}) \rightarrow 0$ and $m^{(2)}_{h}(\mathbf{x}) \rightarrow \pi^{d/2}/2$. Consequently, $A_h(\mathbf{x}) \rightarrow -\pi^{d}$ and $B_h(\mathbf{x}) \rightarrow 0$ in Eq.~\cref{eq:rho_h_quadratic}, which yields $\rho_h(\mathbf{x})^4q(\mathbf{x})^2 \rightarrow 1 + \mathcal{O}(h^2)$.

In addition to the cases of closed manifolds and the interior of a manifold with a boundary, a scenario where a second-order approximation of $\rho_h(\mathbf{x})^2q(\mathbf{x})$ can be easily described arises when $H(\mathbf{x}) = 0$, as with a flat boundary, and $\partial_{\boldsymbol{\eta}_{\mathbf{x}}}q(\mathbf{x}) = 0$, as with a uniform sampling density.
\begin{corollary}
\label{corollary:rho_h_H_zero}
If $\partial \mathcal{M} \neq \emptyset$ and $H(\mathbf{x}) = 0$, then
\begin{equation}
    \rho_h(\mathbf{x})^2q(\mathbf{x}) = C_h(\mathbf{x})(1 + \mathcal{O}(h^2)) \text{ where } C_h(\mathbf{x}) = \frac{B_h(\mathbf{x}) - \sqrt{B_h(\mathbf{x})^2 - 4\pi^dA_{h}(\mathbf{x})}}{2A_h(\mathbf{x})} > 0.
\end{equation}
Moreover, if $\partial_{\boldsymbol{\eta}_{\mathbf{x}}}q(\mathbf{x}) = 0$ then $C_h(\mathbf{x})$ only depends on $b_{\mathbf{x}}$, the distance of $\mathbf{x}$ from the boundary.
\end{corollary}

In general, it seems challenging to obtain a second-order approximation of $\rho_h(\mathbf{x})^2q(\mathbf{x})$ that is independent of $H(\mathbf{x})$ and $\partial_{\boldsymbol{\eta}_{\mathbf{x}}}q(\mathbf{x})$, which are typically unknown for a given dataset. Even when $\mathbf{x}$ lies on the boundary, the terms $A_{h}(\mathbf{x})$ and $B_h(\mathbf{x})$ depend on these quantities as shown in the corollary below.
\begin{corollary}
\label{corollary:rho_h_boundary}
If $\partial\mathcal{M} \neq \emptyset$ then for $\mathbf{x} \in \partial \mathcal{M}$,
\begin{align}
    A_h(\mathbf{x}) &= \frac{\pi^{d-1}}{2}\left(1 - \frac{\pi}{2}\right) + \frac{h(d-1)\pi^{d-1/2}H(\mathbf{x})}{4}\\
    B_h(\mathbf{x}) &= -\frac{\pi^{d-1/2}}{2}h\kappa_{+}(\mathbf{x}).
\end{align}
\end{corollary}
Therefore, for practical purposes, we relax the result in~\Cref{thm:rho_hx_char} to obtain a first-order characterization of $\rho_h(\mathbf{x})^2q(\mathbf{x})$ that depends exclusively on $b_{\mathbf{x}}$, removing the reliance on $H(\mathbf{x})$ and $\partial_{\boldsymbol{\eta}_{\mathbf{x}}}q(\mathbf{x})$, which are typically not accessible for a given dataset.
\begin{corollary}
\label{corollary:rho_h_first_order_char}
If $\partial \mathcal{M} \neq \emptyset$ then
\begin{equation}
    \rho_h(\mathbf{x})^2q(\mathbf{x}) = \zeta_h(\mathbf{x}) +\mathcal{O}(h) \label{eq:rho_h_first_order_char}
\end{equation}
where
\begin{equation}
    \zeta_h(\mathbf{x}) = \frac{\pi^{d/2}}{2}\frac{\frac{b_{\mathbf{x}}}{h}m^{(1)}_{h}(\mathbf{x}) + \sqrt{\left(\frac{b_{\mathbf{x}}}{h}m^{(1)}_{h}(\mathbf{x}) + m^{(0)}_{h}(\mathbf{x})\right)^2 - 2m^{(1)}_{h}(\mathbf{x})^2}}{m^{(2)}_{h}(\mathbf{x})m^{(0)}_{h}(\mathbf{x})-m^{(1)}_{h}(\mathbf{x})^2} \label{eq:zeta_h}
\end{equation}
is a strictly decreasing convex function that depends only on $b_{\mathbf{x}}$, the distance of $\mathbf{x}$ from the boundary, with a range of $\left(1,2\sqrt{\frac{\pi}{\pi-2}}\right]$.
\end{corollary}
The above characterization of the scaling factors $\rho_h(\mathbf{x})$ is utilized in the convergence analysis of our doubly stochastic kernel-based BDE, which is presented in the next section.

%% file: sections/algo.tex
In this section, we first review the standard Gaussian kernel-based BDE proposed in~\cite{berry2017density} and how it is utilized to estimate the boundary points. Next, we propose a BDE that leverages the doubly stochastic kernel obtained via Sinkhorn iterations~\cite{landa2023robust,cuturi2013sinkhorn} applied to the Gaussian kernel, and local PCA. Using~\Cref{thm:rho_hx_char} and~\Cref{corollary:rho_h_first_order_char} we establish the convergence of our BDE to the normal direction to the boundary, up to a factor that depends only on the distance to the boundary. Building on this result, we introduce a straightforward algorithm for computing boundary points of the given data manifold. Finally, we propose a two-parameter family of KDEs that utilize the doubly stochastic kernel and converge to the sampling density up to a scaling factor that depends on the distances to the boundary.

\subsection{Boundary estimation using standard Gaussian kernel}
\label{subsec:std_kernel}
\begin{figure}[h]
    \centering
    \includegraphics[width=0.25\textwidth]{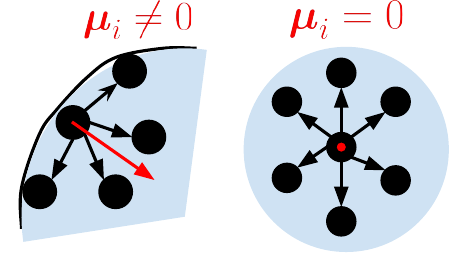}
    \caption{A schematic illustrating $\boldsymbol{\mu}_i$ at the boundary and in the interior.}
    \label{fig:bde_schematic}
\end{figure}

Using the standard Gaussian kernel $\mathbf{K}_{ij}$ as defined in Eq.~\cref{eq:Kij}, Berry and Sauer~\cite{berry2017density} defined a KDE $f$ and a BDE $ \boldsymbol{\mu}$ as follows,
\begin{align}
    f_i &\coloneqq \frac{1}{n-1}\sum_{j=1}^{n}\mathbf{K}_{ij} \label{eq:f}\\
    \boldsymbol{\mu}_i &\coloneqq \frac{1}{n-1}\sum_{j=1}^{n}\mathbf{K}_{ij} (\widetilde{\mathbf{x}}_j-\widetilde{\mathbf{x}}_i).\label{eq:mu}
\end{align}
Intuitively, $\boldsymbol{\mu}_i$ vanishes in the interior while, near the boundary, it is non-zero and points in the direction of $-\boldsymbol{\eta}_{\mathbf{x}_i}$, i.e. opposite to the normal direction to the boundary, as illustrated in~\Cref{fig:bde_schematic}.

Berry and Sauer~\cite{berry2017density} showed in the noiseless setting that $\boldsymbol{\mu}_i$ converges to the normal direction $-\boldsymbol{\eta}_{\mathbf{x}_i}$ up to a global factor that depends on the distance $b_{\mathbf{x}_i}$ of $\mathbf{x}_i$ from the boundary and the sampling density $q(\mathbf{x}_i)$ at $\mathbf{x}_i$. They also showed that $f_i$ converges to the sampling density $q(\mathbf{x}_i)$ up to another global factor that depends on $b_{\mathbf{x}_i}$. Precisely,
\begin{align}
    \mathbb{E}[f_i] &= h^d m_{h}^{(0)}(\mathbf{x}_i)q(\mathbf{x}_i)(1+\mathcal{O}(h)) \text{ and }\\
    \mathbb{E}[\boldsymbol{\mu}_i] &= h^{d+1}\boldsymbol{\eta}_{\mathbf{x}_i} q(\mathbf{x}_i)m_{h}^{(1)}(\mathbf{x}_i) (1+\mathcal{O}(h)).\label{eq:E[mu]}
\end{align}
Taking the ratio of $\mathbb{E}[f_i]$ and $\left\|\mathbb{E}[\boldsymbol{\mu}_i]\right\|$ and thereby eliminating $q(\mathbf{x}_i)$, they obtain
\begin{align}
    &\frac{h\mathbb{E}[f_i]}{\sqrt{\pi}\left\|\mathbb{E}[\boldsymbol{\mu}_i]\right\|}  = \frac{m_{h}^{(0)}(\mathbf{x}_i)(1+\mathcal{O}(h))}{\sqrt{\pi} |m_{h}^{(1)}(\mathbf{x}_i)|} = \left(1+\mathrm{erf}\left(\frac{b_{\mathbf{x}_i}}{h}\right)\right) e^{b_{\mathbf{x}_i}^2/h^2}(1+\mathcal{O}(h)).
\end{align}
The distances to the boundary $\widehat{b}_{\mathbf{x}_i}$ are then estimated by finding the root of the equation $F(b_{\mathbf{x}}) = hf_i/\sqrt{\pi}\left\|\boldsymbol{\mu}_i\right\|$ using Newton's method where $F(b_{\mathbf{x}}) = \left(1+\mathrm{erf}\left(\frac{b_{\mathbf{x}}}{h}\right)\right) e^{b_{\mathbf{x}}^2/h^2}$.

\subsection{Boundary estimation using doubly stochastic kernel and local PCA}
\label{subsec:ds_bde}
Here, we propose a BDE by replacing the standard Gaussian kernel $\mathbf{K}$ in Eq.~(\ref{eq:mu}) with the doubly stochastic kernel $\mathbf{W}$ as defined in Eq.~(\ref{eq:W}). Moreover, instead of using the vectors $\widetilde{\mathbf{x}}_j - \widetilde{\mathbf{x}}_i$ between the original (potentially noisy) data points, our method employs vectors between the PCA-projected points, $\widetilde{\mathbf{U}}_i^T(\widetilde{\mathbf{x}}_j - \widetilde{\mathbf{x}}_i)$ where $\widetilde{\mathbf{U}}_i \in \mathbb{R}^{m \times d}$ serves as an approximation of an orthogonal basis of $T_{\mathbf{x}_i}\mathcal{M}$ and is obtained by applying PCA on a neighborhood of $\widetilde{\mathbf{x}}_i$. At a high level, this design leverages the noise-resilient properties of the doubly stochastic kernel~\cite{landa2023robust,cheng2024bi} combined with the denoising capabilities of local PCA~\cite{gong2010locally}. Overall, our BDE at the $i$-th point is given by,
\begin{align}
    \boldsymbol{\nu}_i &\coloneqq \frac{1}{n-1}\sum_{j=1}^{n}\mathbf{W}_{ij} \widetilde{\mathbf{U}}_i^T(\widetilde{\mathbf{x}}_j-\widetilde{\mathbf{x}}_i). \label{eq:nu}
\end{align}

To maintain the focus of this paper on the doubly stochastic scaling and to simplify our analysis, we adopt the following assumptions in addition to (\refassump{assump:rho}{A1}, \refassump{assump:q}{A2'}, \refassump{assump:manifold}{A3}-\refassump{assump:dim}{A5}):
\begin{enumerate}
    \item[(A6)]\label{assump:tangent} The projector onto the tangent space is assumed to be available. Specifically, we assume access to the matrices $\{\boldsymbol{\mathcal{U}}_i\}_1^n \subset \mathbb{R}^{m \times d}$, where the columns of $\boldsymbol{\mathcal{U}}_i$ form an orthonormal basis of $T_{\mathbf{x}_i}\mathcal{M}$. In our experiments, however, we approximate an orthonormal basis of $T_{\mathbf{x}_i}\mathcal{M}$ using $\widetilde{\mathbf{U}}_i$, obtained by applying local PCA to noisy observations. For readers interested in the analysis showing the convergence of $\widetilde{\mathbf{U}}_i$ to a true basis $\boldsymbol{\mathcal{U}}_i\mathbf{S}_i$ for some $\mathbf{S}_i \in \mathbb{O}(d)$, we refer to the following papers~\cite{singer2012vector,tyagi2013tangent,cheng2013local,kaslovsky2014non,yu2015useful,mohammed2017manifold,aamari2019nonasymptotic}.
    \item[(A7)]\label{assump:ortho_noise} The noise $\boldsymbol{\varepsilon}(\mathbf{x}_i)$ is orthogonal to the tangent space, i.e. $\boldsymbol{\mathcal{U}}_i^T\boldsymbol{\varepsilon}(\mathbf{x}_i) = 0$.
\end{enumerate}

With assumptions \refassump{assump:tangent}{(A6)} and \refassump{assump:ortho_noise}{(A7)}, we are now ready to state our main result, which characterizes our BDE $\boldsymbol{\nu}_i$. The following theorem establishes that our BDE using a doubly stochastic kernel and local PCA converges to the true normal direction $\boldsymbol{\eta}_{\mathbf{x}_i}$ up to a scaling factor that depends only on the distance from the boundary.
\begin{theorem}
\label{theorem:Enu}
Under assumptions (\refassump{assump:rho}{A1}, \refassump{assump:q}{A2'}, \refassump{assump:manifold}{A3}-\refassump{assump:ortho_noise}{A7}), we have
\begin{equation}
    \boldsymbol{\nu}_i = h\beta_h(\mathbf{x}_i) \boldsymbol{\eta}_{\mathbf{x}_i}  + \mathcal{O}\left(h^2\right) + \mathcal{O}^{(h)}(g(m,n))\label{eq:E[nu]}
\end{equation}
where $\mathcal{O}^{(h)}(g(m,n))$ is defined in~\Cref{lem:setup} and
\begin{equation}
    \beta_h(\mathbf{x}_i) = \frac{m_h^{(1)}(\mathbf{x}_i)\zeta_h(\mathbf{x}_i)}{\pi^{d/2} + 2m_h^{(2)}(\mathbf{x}_i)\zeta_h(\mathbf{x}_i)}.
\end{equation}
Here, $-\beta_h(\mathbf{x})$ is a positive and strictly decreasing convex function that depends only on $b_{\mathbf{x}}$, the distance of $\mathbf{x}$ from the boundary, and decays from $\frac{1}{2}(\sqrt{\pi} - \sqrt{\pi-2})$ to $0$ as $b_{\mathbf{x}} \rightarrow \infty$.
\end{theorem}
If the sampling density is uniform, then the convergence order improves from $\mathcal{O}\left(h^2\right)$ to $\mathcal{O}\left(h^3\right)$.

\begin{corollary}
\label{corollary:Enu}
Under assumptions (\refassump{assump:rho}{A1}, \refassump{assump:q}{A2'}, \refassump{assump:tangent}{A3}-\refassump{assump:ortho_noise}{A7}), if $H(\mathbf{x}_i) = 0$ and $q$ is uniform, then
$$\boldsymbol{\nu}_i = h\beta_h(\mathbf{x}_i) \boldsymbol{\eta}_{\mathbf{x}_i}  + \mathcal{O}\left(h^3\right) + \mathcal{O}^{(h)}(g(m,n)).$$
\end{corollary}
The proof follows directly from the proof of~\Cref{theorem:Enu} by substituting zeros for $H$ and the partial derivatives of $q$. 

From the above results, it follows that 
\begin{equation}
    \left\|\boldsymbol{\nu}_i\right\| \approx h|\beta_h(\mathbf{x}_i)| = -h\beta_h(\mathbf{x}_i)
\end{equation}
for sufficiently small $h$, with high probability in the large sample and high-dimensional setting. Since $\beta_h(\mathbf{x})$ depends solely on $b_{\mathbf{x}}$ and is strictly decreasing, we can identify the boundary points by simply thresholding $\{\left\|\boldsymbol{\nu}_i\right\|\}_1^n$, for example by using a chosen percentile of these values. The steps to compute $(\boldsymbol{\nu}_i)_1^n$ are summarized in~\Cref{algo:bde}, and it differs from the Gaussian kernel-based approach~\cite{berry2017density}  described in \Cref{subsec:std_kernel} in two ways: (i) it employs Sinkhorn iterations~\cite{cuturi2013sinkhorn,landa2023robust} to construct the doubly stochastic kernel, and (ii) it eliminates the need for KDE computations and Newton's method to approximate boundary points.
In summary, normalizing the kernel via doubly stochastic scaling reduces the impact of noise~\cite{landa2021doubly,cheng2024bi}, while local PCA provides additional noise filtering, improving the robustness of the estimates.
Our method also maintains computational efficiency by identifying boundary points through a straightforward thresholding of the estimated boundary direction's norm.

\begin{algorithm}[H]
\caption{Boundary direction and points estimator\label{algo:bde}}
\begin{algorithmic}[1]
\REQUIRE $(\widetilde{\mathbf{x}}_i)_1^n \subset \mathbb{R}^m$, $d$ $h$, $k_{\text{nn}}$, $p$.
\STATE Construct $\mathbf{K}$ as in Eq.~(\ref{eq:Kij}).
\STATE Construct $\mathbf{W}$ as in Eq.~(\ref{eq:W}) using Sinkhorn iterations~\cite{cuturi2013sinkhorn,landa2023robust}.
\FOR{$i = 1 \rightarrow n$}
    \STATE $\widetilde{\mathcal{N}}_{i} \leftarrow k_{\text{nn}}$-nearest neighbors of $\widetilde{\mathbf{x}}_i$.
    \STATE $\widetilde{\mathbf{U}}_i \leftarrow$ PCA on $\{\widetilde{\mathbf{x}}_j: j \in \widetilde{\mathcal{N}}_{i}\}$.
    \STATE Set $\boldsymbol{\nu}_i = \frac{1}{n-1}\sum_{j = 1}^{n}\mathbf{W}_{ij}\widetilde{\mathbf{U}}_i^T(\widetilde{\mathbf{x}}_j-\widetilde{\mathbf{x}}_i)$.
\ENDFOR
\STATE $\widehat{\mathcal{B}}_{p} = \{k: \left\|\boldsymbol{\nu}_k\right\| \geq (100-p)\text{th percentile of } (\left\|\boldsymbol{\nu}_i\right\|)_1^n\}$.
\end{algorithmic}
\end{algorithm}

Having detected the boundary points, we can easily estimate distances to the boundary by computing the minimum of the shortest path distances of each point from the boundary points. To ensure robustness to noise, instead of using naive Euclidean distances between nearest neighbors as input to the shortest path algorithm, we use Euclidean distances computed from the PCA-projected embeddings of local neighborhoods. The pseudocode is provided below.

\begin{algorithm}[h]
\caption{Distance to boundary estimator\label{algo:dist_to_bdry}}
\begin{algorithmic}[1]
\REQUIRE $(\widetilde{\mathbf{x}}_i)_1^n \subset \mathbb{R}^m$, $d$, $k_{\text{nn}}$, $\widehat{\mathcal{B}}_{p}$.
\FOR{$i = 1 \rightarrow n$}
    \STATE $\widetilde{\mathcal{N}}_{i} \leftarrow k_{\text{nn}}$-nearest neighbors of $\widetilde{\mathbf{x}}_i$.
    \STATE $\widetilde{\mathbf{U}}_i \leftarrow$ PCA on $\{\widetilde{\mathbf{x}}_j: j \in \widetilde{\mathcal{N}}_{i}\}$.
    \STATE $\mathrm{dist}_{ij} \leftarrow \left\|\widetilde{\mathbf{U}}_i^T(\widetilde{\mathbf{x}}_j-\widetilde{\mathbf{x}}_i)\right\|_2$ for all $j \in \widetilde{\mathcal{N}}_{i}$.
\ENDFOR
\STATE $\mathrm{shortest\_dist}_{ij} \leftarrow \mathrm{Dijkstra}(\mathrm{dist}, \mathrm{directed = False})$ for all $i \in [1,n]$ and $j \in \widehat{\mathcal{B}}_{p}$.
\STATE $\widehat{b}_{\mathbf{x}_i} \leftarrow \min_{j\in \mathcal{S}_{\partial \mathcal{M}}}\mathrm{shortest\_dist}_{ij}$.
\end{algorithmic}
\end{algorithm}

%% file: sections/kde.tex
Standard kernel density estimators for distributions with bounded support are well known to exhibit bias near the boundary \cite{jones1996simple,marron1994transformations,cowling1996pseudodata,karunamuni2005boundary}. In particular, without explicit correction, kernel mass extends beyond the boundary of the support, leading to systematic underestimation of the density in regions close to the boundary. This phenomenon also persists in the manifold setting. However, having identified the boundary in the previous section, we can utilize the distance to the boundary to construct an explicit correction. Here, we introduce a two-parameter family of kernel density estimators that extends the method of \cite{landa2023robust}, originally developed for closed manifolds, by incorporating a boundary correction term based on $m^{(0)}_{h/\sqrt{s}}(\mathbf{x})$ and $\zeta_h(\mathbf{x})$, as defined in Eq.~(\ref{eq:m0}, \ref{eq:zeta_h}). The resulting estimator is given by:
\begin{equation}
    \widehat{q}_i \coloneqq \frac{1}{(n-1)h^{d}}\left(\frac{1}{\pi^{ds/2}s^{d/2}}\frac{1}{(m^{(0)}_{h/\sqrt{s}}(\mathbf{x}_i)\zeta_{h}(\mathbf{x}_i)^{s})^u}\sum_{j=1}^{n}\frac{\mathbf{W}_{ij}^s}{(m^{(0)}_{h/\sqrt{s}}(\mathbf{x}_j)\zeta_{h}(\mathbf{x}_j)^{s})^{1-u}}\right)^{1/(1-s)}
    \label{eq:kde}
\end{equation}
where $s > 0, s \neq 1$ and $u \in (0,1)$.
In the absence of the terms involving $m^{(0)}_{h/\sqrt{s}}$ and $\zeta_h$, the above formulation reduces to the estimator in~\cite{landa2023robust} for closed manifolds (up to a global constant that depends on $h$, $s$ and $d$). The following theorem establishes the convergence of the above KDE to the true sampling density.
\begin{theorem}
\label{theorem:Eq}
    Under assumptions (\refassump{assump:rho}{A1}, \refassump{assump:q}{A2'}, \refassump{assump:manifold}{A3}-\refassump{assump:dim}{A5}), we have
    \begin{equation}
        \widehat{q}_i   = \left(q(\mathbf{x}_i) + \mathcal{O}(h)\right)(1 + \mathcal{O}^{(h)}(g(m,n))),
    \end{equation}
    where $\mathcal{O}^{(h)}(g(m,n))$ is defined in~\Cref{lem:setup}.
\end{theorem}
It is worth noting that for closed manifolds with intrinsic dimension $d \leq 3$, the estimator in \cite{landa2023robust} achieves a second-order error of $\mathcal{O}(h^2)$\textemdash a result conjectured to hold for higher dimensions as well. In contrast, our estimator for manifolds with boundary is limited to first-order accuracy $\mathcal{O}(h)$. A primary cause of this limitation is that the mean curvature $H(\mathbf{x})$ induced by the boundary is unknown (see~\cref{thm:rho_hx_char}), thereby preventing the implementation of higher-order corrections. Nevertheless, the proposed first-order estimator in Eq.~(\ref{eq:kde}) is computationally feasible, and is implemented using the estimated distances to the boundary from~\Cref{algo:dist_to_bdry} to compute the correction terms in $\hat{q}_i$.


%% file: sections/exp.tex
In \Cref{sec:robust_boundary_estimation}, we estimate points located on and near the boundary for data uniformly sampled from two geometries: (a) a circular annulus and (b) a curved truncated torus. These estimations are conducted under three noise scenarios: no noise, homoskedastic noise, and heteroskedastic noise.
In \Cref{subsec:kde-exp}, we demonstrate how our boundary correction term improves density estimation on the circular annulus and truncated torus under heteroskedastic noise.
Finally, in~\Cref{subsec:mnist}, we apply our method to high-dimensional MNIST dataset to identify images that lie close to and far from the boundary. In the following experiments, we evaluate the performance of the following methods:
\begin{enumerate}[leftmargin=*]
    \item The standard Gaussian heat kernel as in~\Cref{subsec:std_kernel}, and a local PCA enhanced variant. In the latter, during the calculation of $\boldsymbol{\mu}_i$ in Eq.~\cref{eq:mu}, $\widetilde{\mathbf{x}}_j-\widetilde{\mathbf{x}}_i$ is replaced with $\widetilde{\mathbf{U}}_i^T(\widetilde{\mathbf{x}}_j-\widetilde{\mathbf{x}}_i)$ as in Eq.~\cref{eq:nu}, by employing local PCA in a small neighborhood around $\widetilde{\mathbf{x}}_i$. These methods are denoted as \textalgo{Gaussian} and \textalgo{Gaussian+LPCA}, respectively.
    \item Our approach using the doubly stochastic kernel with local PCA as in~\Cref{algo:bde}, and a variant of it which, while calculating $\boldsymbol{\nu}_i$ in Eq.~\cref{eq:nu}, directly uses the difference in the observations $\widetilde{\mathbf{x}}_j-\widetilde{\mathbf{x}}_i$ as in Eq.~\cref{eq:mu}, instead of the local PCA projections. These methods are referred to as \textalgo{DS+LPCA} and \textalgo{DS}, respectively.
    \item The method proposed in~\cite{calder2022boundary} that estimates the boundary using an approach similar to \textalgo{Gaussian}, but replaces the Gaussian kernel with a binary cutoff kernel. We refer to this method as~\textalgo{Binary}. The implementation of this method is available in~\cite{graphlearning}.
\end{enumerate}
To estimate the boundary, we threshold the estimated distances $(\widehat{b}_{\mathbf{x}_i})_1^n$ and the norm of BDE $(\left\|\boldsymbol{\nu}_i\right\|)_1^n$ based on a specified percentile. For clarity, we define the ground truth $p$-th percentile boundary as
\begin{equation}
    \mathcal{B}_{p} = \{k: b_{\mathbf{x}_k} \leq p\text{-th percentile of } (b_{\mathbf{x}_i})_1^n\}\label{eq:Cp}.
\end{equation}
For the methods \textalgo{Gaussian}, \textalgo{Gaussian+LPCA} and \textalgo{Binary}, the 
estimated $p$-th percentile boundary is set to 
\begin{equation}
    \widehat{\mathcal{B}}_{p} = 
    \{k: \widehat{b}_{\mathbf{x}_k} \leq p\text{-th percentile of }  (\widehat{b}_{\mathbf{x}_i})_1^n\},
    \label{eq:Chatp2}
\end{equation}
and for \textalgo{DS} and \textalgo{DS+LPCA}, it is set to 
\begin{equation}
    \widehat{\mathcal{B}}_{p} =
    \{k: \left\|\boldsymbol{\nu}_k\right\| \geq (100-p)\text{th percentile of } (\left\|\boldsymbol{\nu}_i\right\|)_1^n\}. 
    \label{eq:Chatp}
\end{equation}
To quantify the discrepancy between the estimated and true boundaries, we use the Jaccard index,
\begin{equation}
    \mathcal{J}_p = \frac{|\mathcal{B}_{p} \cap \widehat{\mathcal{B}}_{p}|}{|\mathcal{B}_{p} \cup \widehat{\mathcal{B}}_{p}|}. \label{eq:jaccard}
\end{equation}
The code is available at \url{https://github.com/chiggum/robust_boundary_estimation}.

\subsection{Robust boundary estimation}
\label{sec:robust_boundary_estimation}
\subsubsection{Data descriptions}
\label{subsec:data_descriptions}
\begin{figure}[h]
    \centering
    \includegraphics[width=0.6\linewidth]{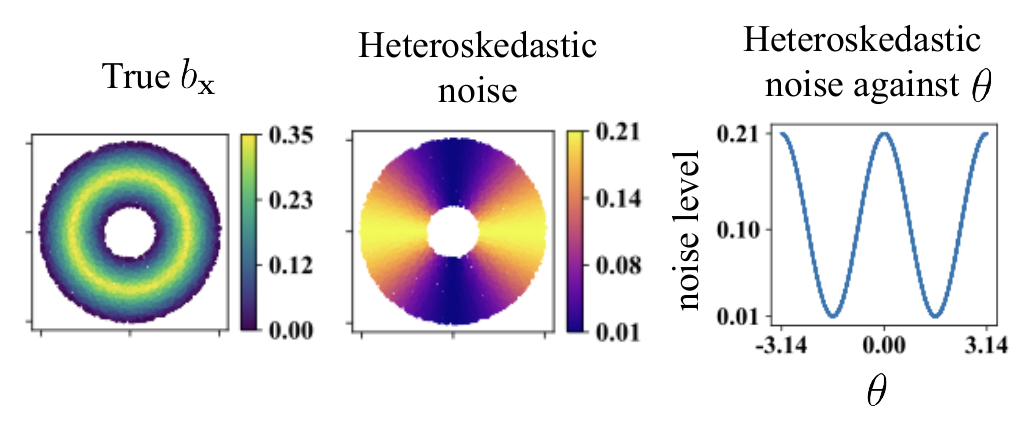}
    \caption{From left to right: (i) the true distances of uniformly sampled points from the boundary of the circular annulus, (ii) the noise levels visualized on the annulus, and (iii) the noise levels plotted as a function of the angular coordinate $\theta$ of the data points.}
    \label{fig:annulus0}
\end{figure}

For our first set of datasets, we generate $n=7076$ uniformly distributed points in $\mathbb{R}^2$, sampled as $\{(r_i\cos(\theta_i),r_i\sin(\theta_i))\}_1^n$ from a circular annulus with an outer radius of $1$ and an inner radius of $0.3$, i.e. $r_i \in [0.3,1]$ and $\theta_i \in (-\pi,\pi]$ (see Figure~\ref{fig:annulus0}). We then project these points into $m=2000$ dimensions using an orthogonal transformation. This forms our noiseless dataset. Let $\mathcal{S}$ be the $2$-dimensional subspace containing the annulus. We construct two noisy variants of the data (a) by adding homoskedastic noise in a ball of radius $0.2$ at each point, ensuring that the noise is orthogonal to $\mathcal{S}$ (this preserves the ground truth distances of the points from the boundary for evaluation) and (b) by adding heteroskedastic noise in a ball of radius $0.01 + 0.2(1+\cos(2\theta_i))/2$ at each point, also orthogonal to $\mathcal{S}$. The noise distribution is illustrated in Figure~\ref{fig:annulus0}.

\begin{figure}[h]
    \centering
    \includegraphics[width=\linewidth]{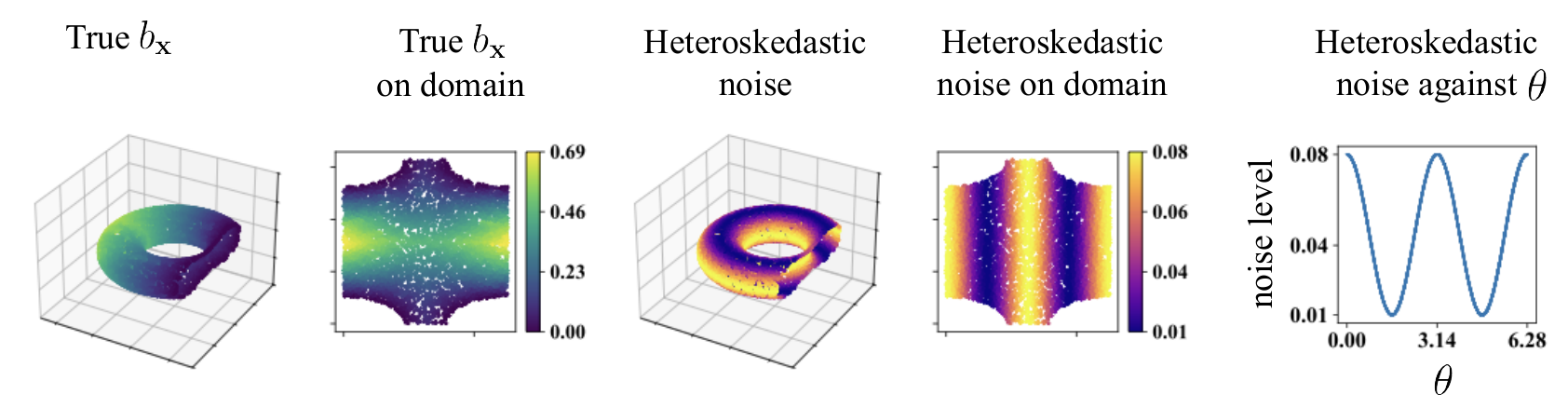}
    \caption{From left to right: (i) true distances of uniformly distributed points from the boundary of a truncated torus, (ii) true distances to the boundary visualized on the two-dimensional domain of the truncated torus, (iii) noise level illustrated on the truncated torus, (iv) noise level visualized on the domain, and (v) noise level as a function of the angular coordinate $\theta$ of the data points.}
    \label{fig:torus0}
\end{figure}

For our second set of datasets, we generate $n=6163$ uniformly distributed points 
$$\{((R + r\cos(\theta_i))\cos(\phi_i),(R + r\cos(\theta_i))\cos(\phi_i),r\sin(\theta_i))\}_1^n$$
on a curved torus in $\mathbb{R}^{3}$ where $R=0.25$, $r=1/4\pi^2R$, $\theta_i \in [0,2\pi)$ and $\phi_i \in [0,2\pi)$. The torus is then truncated by removing points with $x$-coordinate greater than $0.18$, as shown in Figure~\ref{fig:torus0}. These points are projected into $m=2000$ dimensions using an orthogonal transformation, forming our noiseless dataset. Let $\mathcal{S}$ represent the $3$-dimensional subspace containing the truncated torus. We generate two noisy variants of the data (a) by adding homoskedastic noise in a ball of radius $0.075$ at each point while ensuring that the noise is orthogonal to $\mathcal{S}$ and (b) by adding heteroskedastic noise in a ball of radius $0.01 + 0.075(1+\cos(2\theta_i))/2$ at each point, also orthogonal to $\mathcal{S}$. The noise distribution is illustrated in Figure~\ref{fig:torus0}.

\subsubsection{Interpreting results}
\label{subsec:interpret}
For both datasets, we construct the standard Gaussian kernel and the doubly stochastic kernel~\cite{landa2023robust} using a bandwidth of $h=0.1$. Using these kernels, $d=2$ and $k_{\text{nn}}=256$ neighbors for local PCA computations, we estimate the distances to the boundary using \textalgo{DS+LPCA}, \textalgo{DS}, \textalgo{Gaussian+LPCA} and \textalgo{Gaussian}. We use the same number of nearest neighbors while estimating distances to boundary using \textalgo{Binary}.  Subsequently we compute the $p$-th percentile boundary $\widehat{\mathcal{B}}_p$ for each method (Eq.~\cref{eq:Chatp2}-\cref{eq:Chatp}). To compare the performance, we calculate the Jaccard index $\mathcal{J}_p$ (Eq.~\cref{eq:jaccard}) against the ground truth boundary $\mathcal{B}_p$ (Eq.~\cref{eq:Cp}) for various values of $p$. The results are presented in~\Cref{fig:jaccard_annulus,fig:jaccard_torus,fig:bde}, and we summarize key observations from the figures below.

\begin{enumerate}[leftmargin=*]
    \item Noiseless setting: 
    \begin{enumerate}[leftmargin=0.5cm]
        \item Since the circular annulus is flat, local PCA has no impact. As a result, \textalgo{DS+LPCA} and \textalgo{DS} achieve identical Jaccard indices, while \textalgo{Gaussian+LPCA} and \textalgo{Gaussian} also perform equivalently. All the methods seem to achieve similar performance while the \textalgo{Binary} method has slightly better performance. The slightly better performance of \textalgo{DS} over \textalgo{Gaussian} can be attributed to the norm of the BDE, $\left\|\boldsymbol{\nu}_i\right\|$ in Eq.~\cref{eq:nu} and $\left\|\boldsymbol{\mu}_i\right\|$ in Eq.~\cref{eq:mu}. The BDE norm is expected to be zero in the interior and non-zero near the boundary of the manifold (see~\Cref{fig:bde_schematic,fig:bde}). However, \textalgo{Gaussian} frequently yields non-zero BDE, whereas \textalgo{DS} maintains the BDE close to zero within the annulus interior.
        \item The curvature of the torus makes local PCA an important factor. We observe that \textalgo{DS+LPCA} and \textalgo{Gaussian+LPCA} achieve comparable performance and outperform \textalgo{DS} and \textalgo{Gaussian}, indicating that local PCA enhances the accuracy of boundary estimation. At the same time, the performance of \textalgo{Binary} degrades significantly, suggesting the method is sensitive to the curvature. This is consistent with the observation in~\cite{calder2022boundary} that the method \textalgo{Binary} is designed for Euclidean domains, not general smooth manifolds.
    \end{enumerate}
    
    \begin{figure}[t]
        \centering
        \setlength{\tabcolsep}{0pt}
        \begin{tabular}{c}
            \textbf{Circular Annulus}\\
            \includegraphics[width=\textwidth]{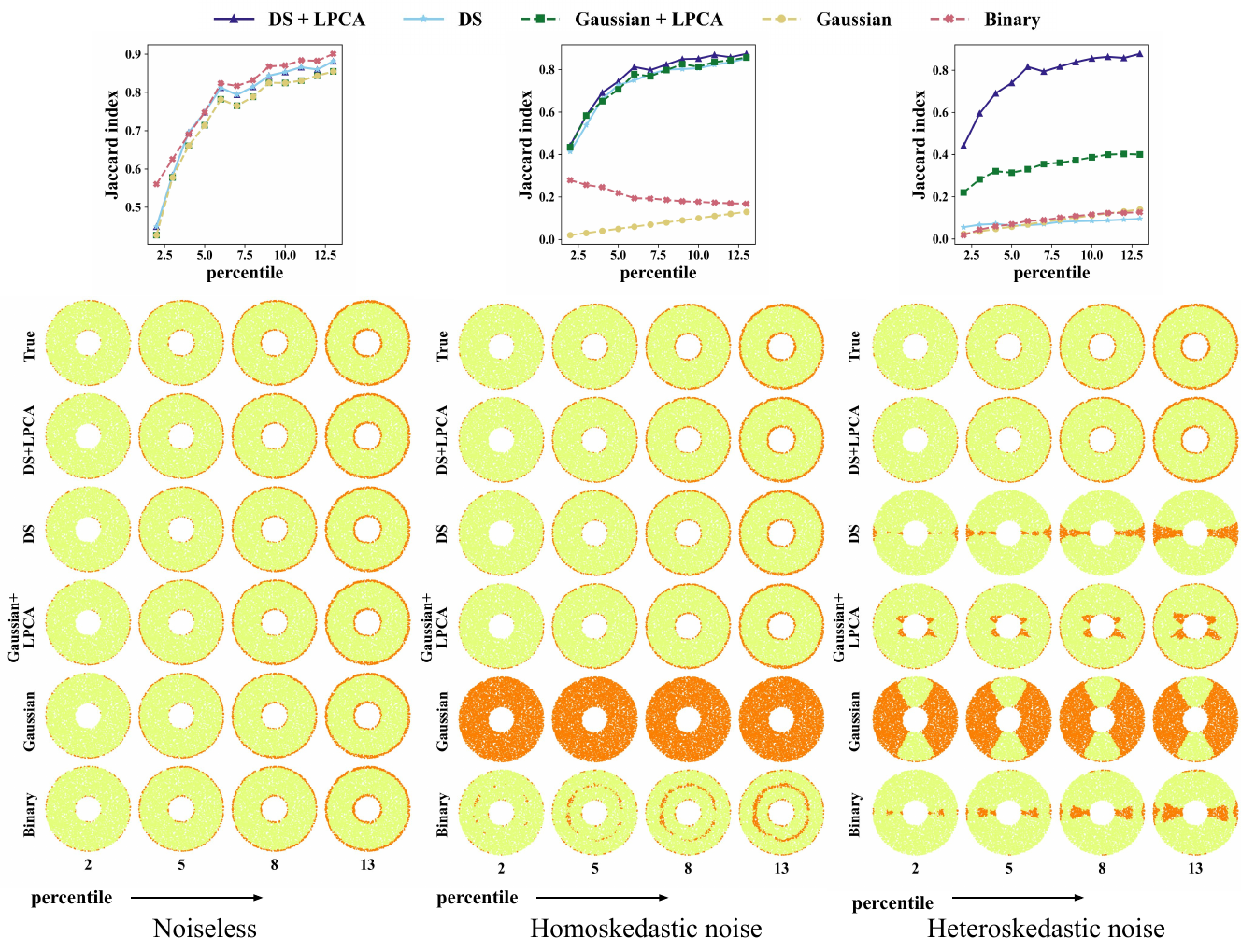}
        \end{tabular}
        \caption{(Top) The Jaccard index (Eq.~\cref{eq:jaccard}) is shown between the ground truth $p$-th percentile boundary $\mathcal{B}_p$ and the estimated $p$-th percentile boundary $\widehat{\mathcal{B}}_p$ for $p \in [2,14]$. (Bottom) Visualizations of the boundaries $\mathcal{B}_p$ and $\widehat{\mathcal{B}}_p$ (in orange) are provided for competing methods at $p \in \{2,5,8,11\}$.}
        \label{fig:jaccard_annulus}
    \end{figure}

    \begin{figure}[t]
    \centering
    \setlength{\tabcolsep}{0pt}
    \begin{tabular}{c}
        \textbf{Truncated torus}\\
        \includegraphics[width=\textwidth]{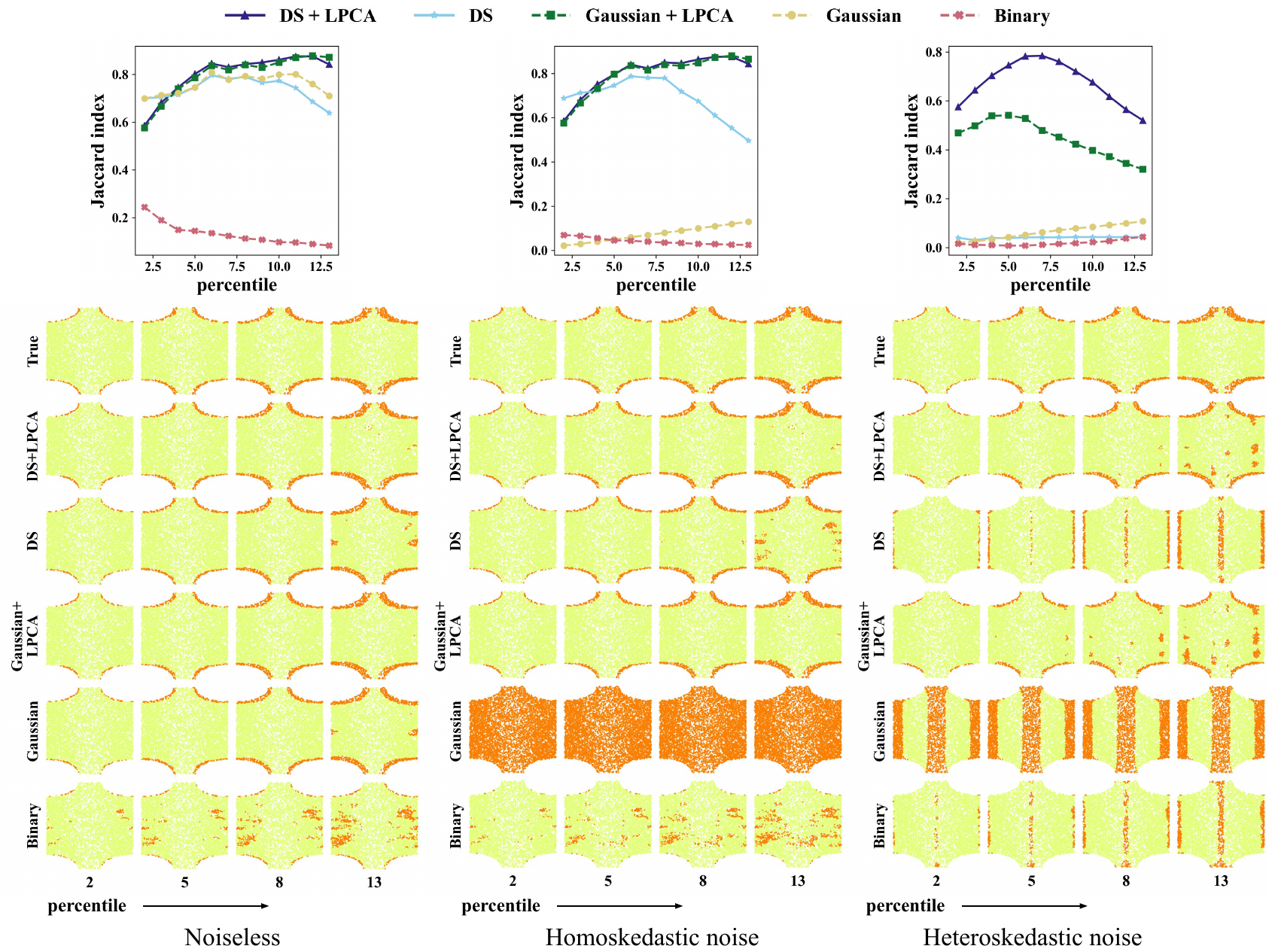}
    \end{tabular}
    \vspace{0.2cm}
    \captionsetup{width=\linewidth}
    \caption{(Top) The Jaccard index (Eq.~\cref{eq:jaccard}) is shown between the ground truth $p$-th percentile boundary $\mathcal{B}_p$ and the estimated $p$-th percentile boundary $\widehat{\mathcal{B}}_p$ for $p \in [2,14]$. (Bottom) Visualizations of the boundaries $\mathcal{B}_p$ and $\widehat{\mathcal{B}}_p$ (in orange) are provided for competing methods at $p \in \{2,5,8,11\}$. The visualizations are presented on the two-dimensional domain of the truncated torus. }
    \label{fig:jaccard_torus}
\end{figure}

    \item Under homoskedastic noise: For both the circular annulus and the truncated torus, the Jaccard index for \textalgo{DS+LPCA}, \textalgo{DS}, and \textalgo{Gaussian+LPCA} remains robust against noise, whereas the performance of \textalgo{Gaussian} and \textalgo{Binary} deteriorates significantly. The degradation of \textalgo{Gaussian} can be attributed to the behavior of the BDE observed in~\Cref{fig:bde} where the BDE from \textalgo{Gaussian} is completely corrupted by noise. Specifically, the BDE is relatively large in the interior and smaller at the boundary, contrary to the expected trend as illustrated in~\Cref{fig:bde_schematic}. In fact, the norm of the BDE is so large that the estimated $\widehat{b}_{\mathbf{x}_i}$ collapse to zero (due to Newton's method converging to negative values that are clipped to zero), resulting in degenerate boundary in~\Cref{fig:jaccard_annulus}. Interestingly, while the norm of the BDE for \textalgo{DS} increases in the interior compared to the noiseless case, the boundary values remain relatively higher. This ensures that the estimates $\left\|\boldsymbol{\nu}_i\right\|$ at the boundary are higher than those in the interior, enabling more accurate estimates of boundary points. Lastly, as shown in~\Cref{fig:bde}, \textalgo{DS+LPCA} and \textalgo{Gaussian+LPCA} exhibit similar BDE values, leading to comparable boundary estimates as shown in~\Cref{fig:jaccard_annulus}, consistent with the noiseless case. This again highlights that incorporating local PCA significantly enhances the accuracy of boundary estimation, as evidenced by the substantial increase in the Jaccard index from \textalgo{Gaussian} to \textalgo{Gaussian+LPCA}.

    \begin{figure}[t]
        \centering
        \setlength{\tabcolsep}{0pt}
        \begin{tabular}{cc}
        \textbf{Circular Annulus} & \textbf{Truncated torus}\\
        \includegraphics[width=0.5\textwidth]{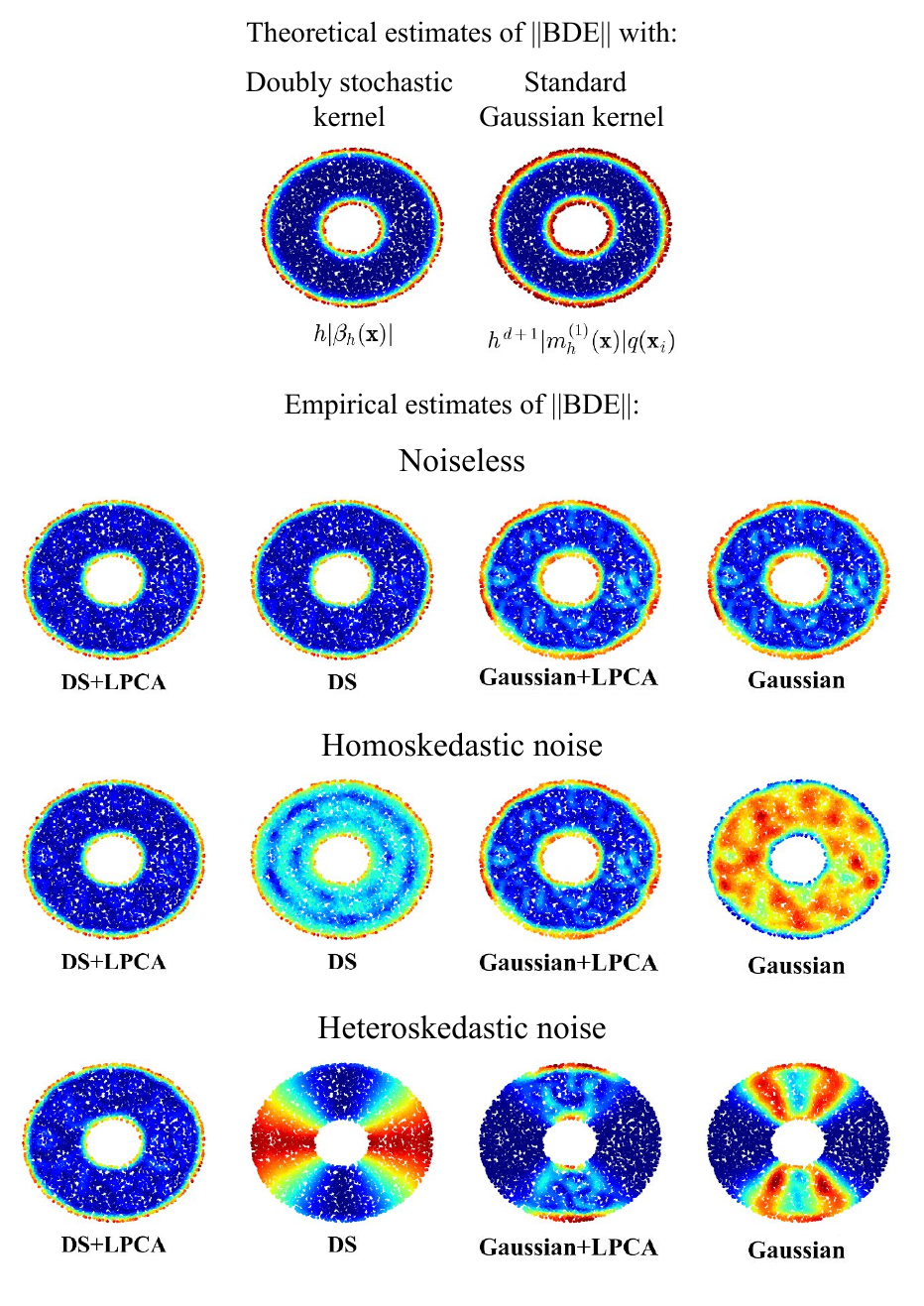} & \includegraphics[width=0.5\textwidth]{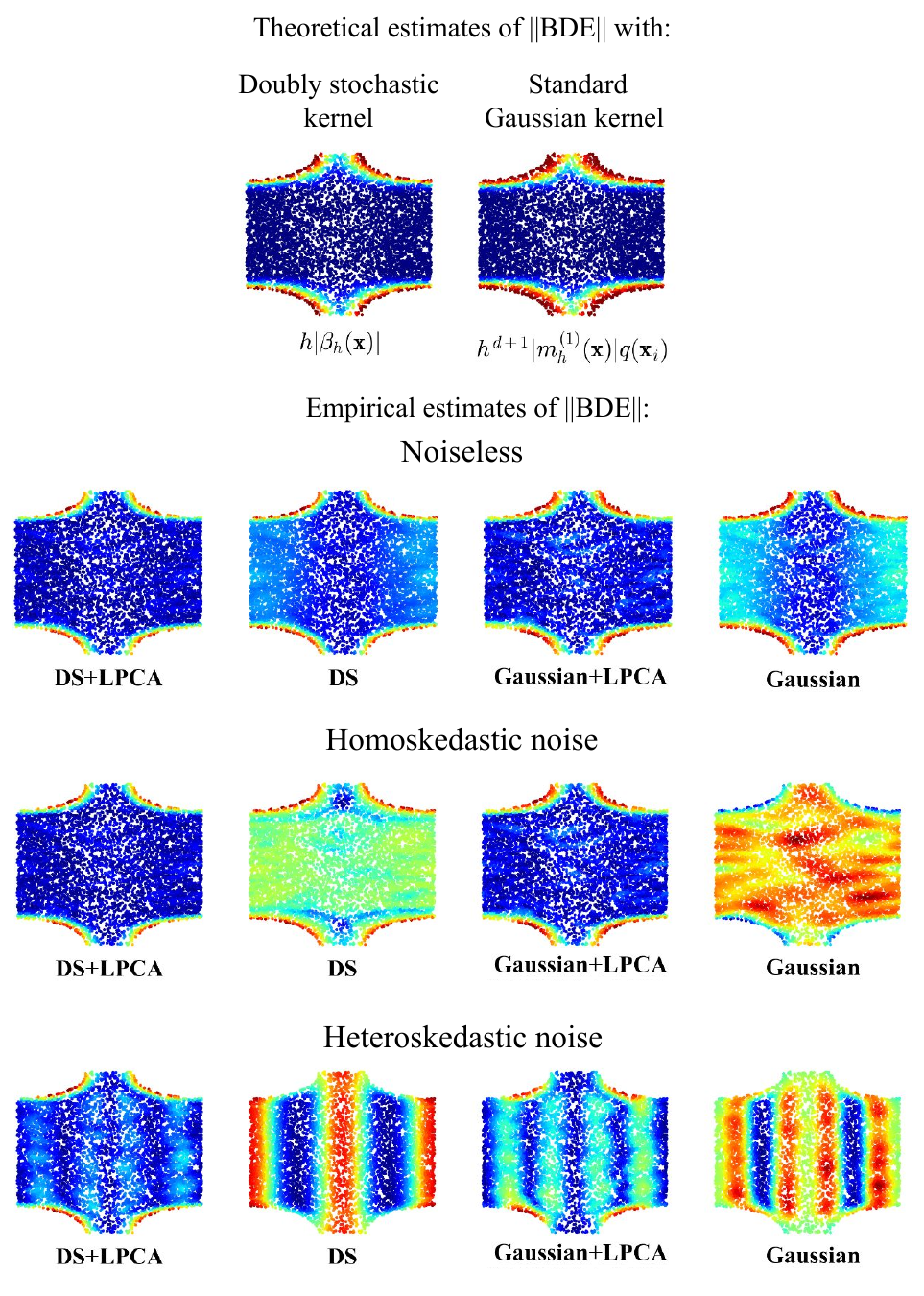}\\
        \end{tabular}
        \caption{(Top) The theoretical norm of the BDE (Eq.~\cref{eq:E[mu]} and Eq.~\cref{eq:E[nu]}) at each point on the data manifold are visualized. These are computed using the true distances to the boundary and uniform density. (Bottom) The empirical norm of the BDE due to each method (Eq.~\cref{eq:mu} and Eq.~\cref{eq:nu}) are visualized on the data manifold. For each case, the values are normalized by the maximum. For the truncated torus, the visualizations are presented onto its two-dimensional domain.}
        \label{fig:bde}
    \end{figure}

    \item Under heteroskedastic noise: For both the circular annulus and the truncated torus, \textalgo{DS+LPCA} yields significantly high Jaccard index compared to \textalgo{Gaussian+LPCA}, while the performance of \textalgo{DS} deteriorates significantly. In fact, \textalgo{Gaussian} and \textalgo{DS} identify the points in the regions with high noise levels (see~\Cref{fig:annulus0,fig:torus0} for noise distributions) as the boundary points as shown in~\Cref{fig:jaccard_annulus,fig:jaccard_torus}.
     This is again attributed to the corruption of BDE due to noise as shown in~\Cref{fig:bde}. Although the BDE for \textalgo{Gaussian+LPCA} is also affected by heteroskedastic noise, the impact is less pronounced compared to \textalgo{DS} and \textalgo{Gaussian}. Nonetheless, the noise still leads to inaccuracies in the boundary estimates as shown in~\Cref{fig:jaccard_annulus,fig:jaccard_torus}. In contrast, the BDE for \textalgo{DS+LPCA} shows greater resilience to noise, indicating that both doubly stochastic scaling and local PCA are crucial for producing robust boundary estimates.
\end{enumerate}

\subsection{Boundary-corrected density estimation}
\label{subsec:kde-exp}

\subsubsection{Datasets}
\label{subsubsec:kde-data}

Once again, we sampled points $\{(r_i\cos(\theta_i),r_i\sin(\theta_i))\}_1^n$ from a circular annulus with $r_i \in [0.3,1]$ and $\theta_i \in (0,2\pi]$. This time, however, we sample $n = 7076$ points with non-uniform density where $\theta_i$ is sampled from $\text{Normal}(0, 0.25\pi^2)$ modulo $2\pi$ and $r_i \sim \text{Uniform}(0.3,1)$. We use the heteroskedastic noise conditions $0.01 + 0.2(1+\cos(2\theta_i))/2$ from \Cref{subsec:data_descriptions}. For the second dataset, we sample $n = 11,175$ points from a truncated curved torus with heteroskedastic noise in a ball of radius $0.01 + 0.075(1+\cos(2\theta_i))/2$, this time with a non-uniform density where $\theta_i \sim \text{Truncated-Normal}((-\pi,\pi),0,0.75\pi)$ and $\phi_i \sim \text{Uniform}(0,2\pi)$.

\subsubsection{Results}
\label{subsubsec:kde-results}
\begin{figure}[h!]
    \centering
    \includegraphics[width=0.99\textwidth]{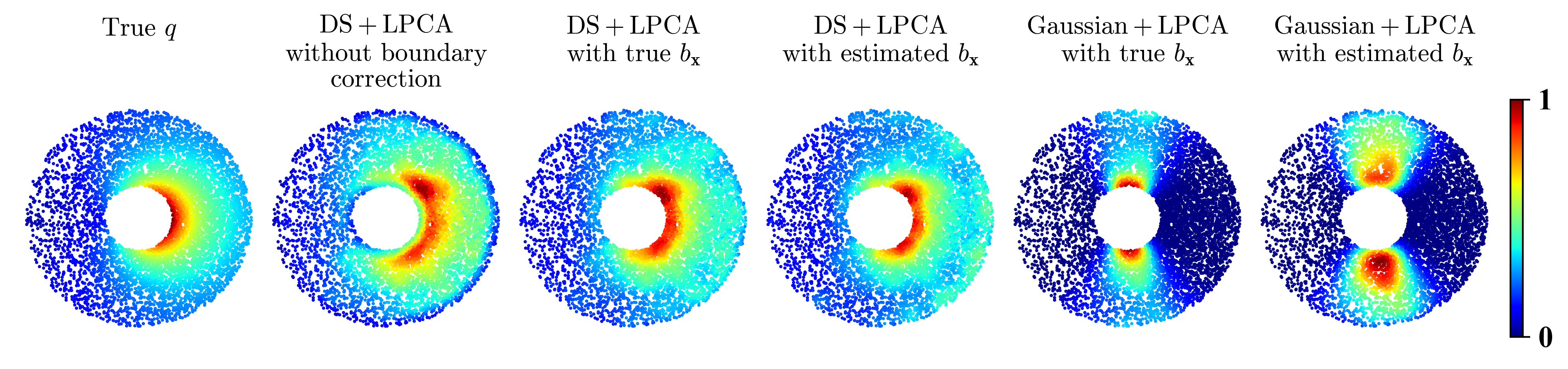}
    \caption{Comparison of kernel density estimators on a circular annulus with heteroskedastic noise compared to the ground truth density $q$ (left). Each estimate is normalized by its respective maximum value.}
    \label{fig:kde_annulus}
\end{figure}

\begin{figure}[h!]
    \centering
    \includegraphics[width=0.5\textwidth]{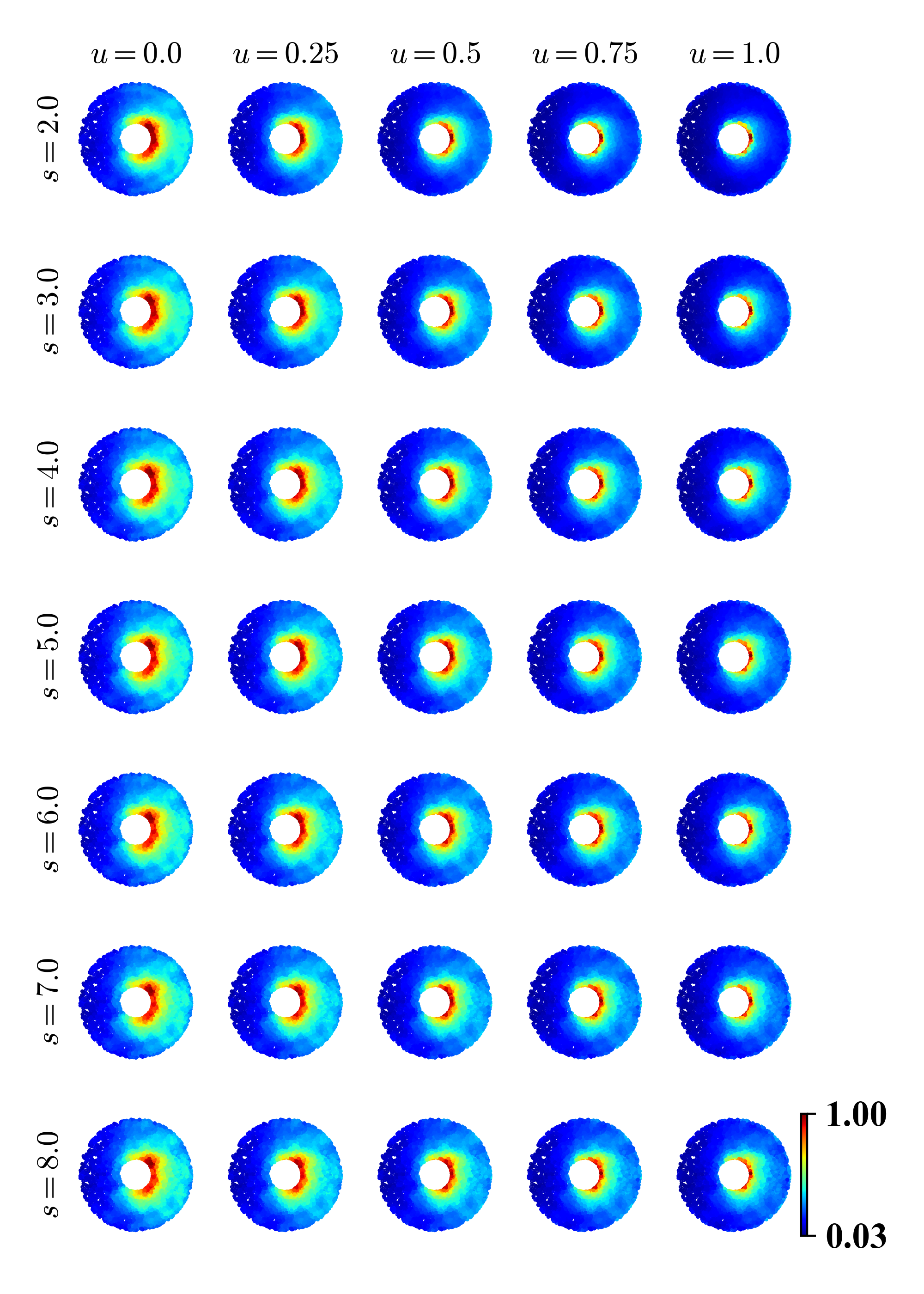}\includegraphics[width=0.5\textwidth]{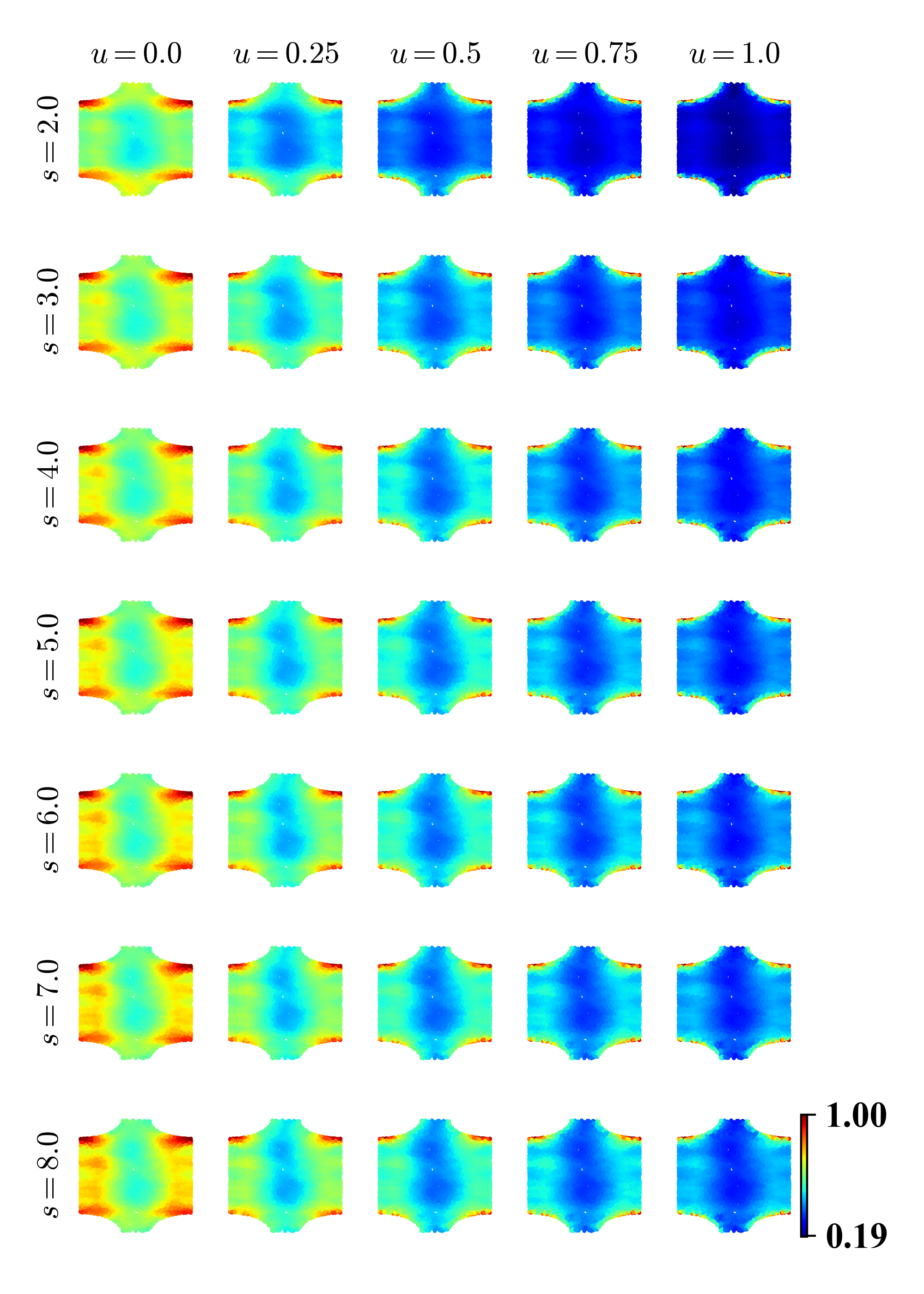}
    \caption{The two-parameter family of kernel density estimators defined in Eq.~(\ref{eq:kde}), evaluated for various values of $s$ and $u$ on a circular annulus and a torus with heteroskedastic noise, in the setting where the true distances to the boundary are known \emph{a priori}. As before, each estimate is normalized by its respective maximum value.}
    \label{fig:kde_sweep}
\end{figure}

For the annulus, we set the intrinsic dimension to $d = 2$, the bandwidth to $h = 0.1$, and the number of nearest neighbors to $k_{\mathrm{nn}} = 1024$ for both the doubly stochastic and standard Gaussian kernels. The same kernels are used for both boundary estimation and kernel density estimation. For the KDE defined in Eq.~(\ref{eq:kde}), we fix the parameters to $s = 8.0$ and $u = 0$; an ablation study over these parameters is provided in \Cref{fig:kde_sweep}. As shown in \Cref{fig:kde_annulus}, the KDE without boundary correction systematically underestimates the density near the boundary due to kernel mass extending beyond the support, as discussed in \Cref{subsec:ds_kde}. The boundary correction term in Eq.~(\ref{eq:kde}) effectively resolves this bias. In contrast, the standard Gaussian kernel–based KDE with boundary correction \cite{berry2017density} remains corrupted by heteroskedastic noise.
\begin{figure}[h!]
    \centering
    \includegraphics[width=0.99\textwidth]{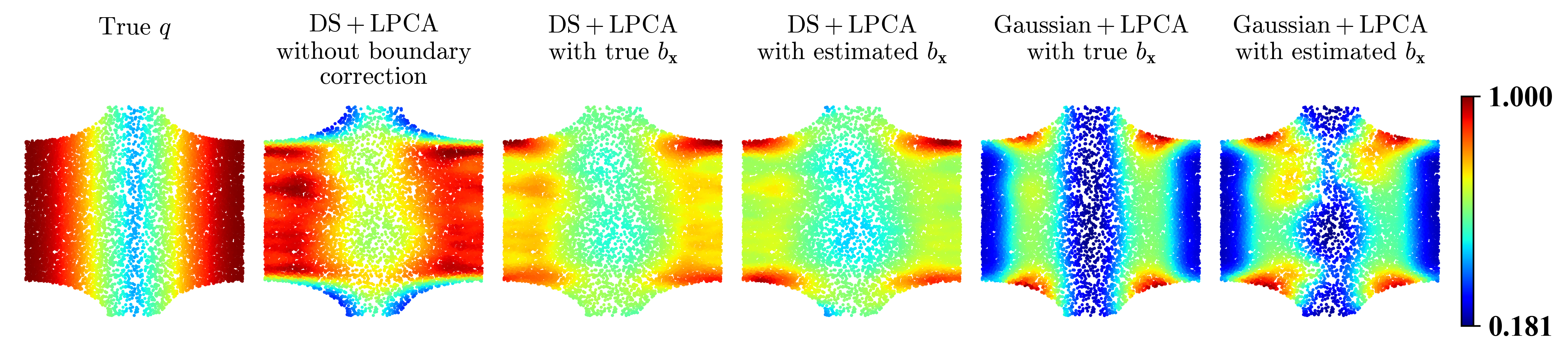}
    \caption{Comparison of kernel density estimators on truncated torus with heteroskedastic noise compared to the ground truth density $q$ (left). Each estimate is normalized by its respective maximum value.}
    \label{fig:kde_torus}
\end{figure}

For the truncated torus, we again set $d = 2$ and $h = 0.1$. Boundary estimation is performed using $k_{\mathrm{nn}} = 256$, after which the KDE is computed with $k_{\mathrm{nn}} = 1024$ once the distances to the boundary have been estimated. As in the annulus case, we fix $s = 8$ and $u = 0$ in Eq.~(\ref{eq:kde}), with a corresponding ablation study reported in \Cref{fig:kde_sweep}. The results in \Cref{fig:kde_torus} again demonstrate that omitting boundary correction leads to an underestimation of the density near the boundary. In contrast to the annulus, the proposed boundary-corrected KDE in Eq.~(\ref{eq:kde}) exhibits a mild overestimation near the boundary; nevertheless, it consistently outperforms the standard Gaussian kernel–based KDE with boundary correction.

\subsection{MNIST} 
\label{subsec:mnist}
\begin{figure}[h!tbp]
    \centering
    \begin{tabular}{c}
       \includegraphics[width=0.9\linewidth]{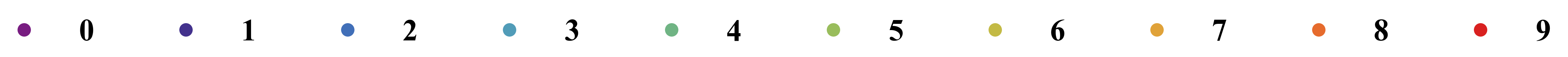}\\
        \includegraphics[width=0.95\linewidth]{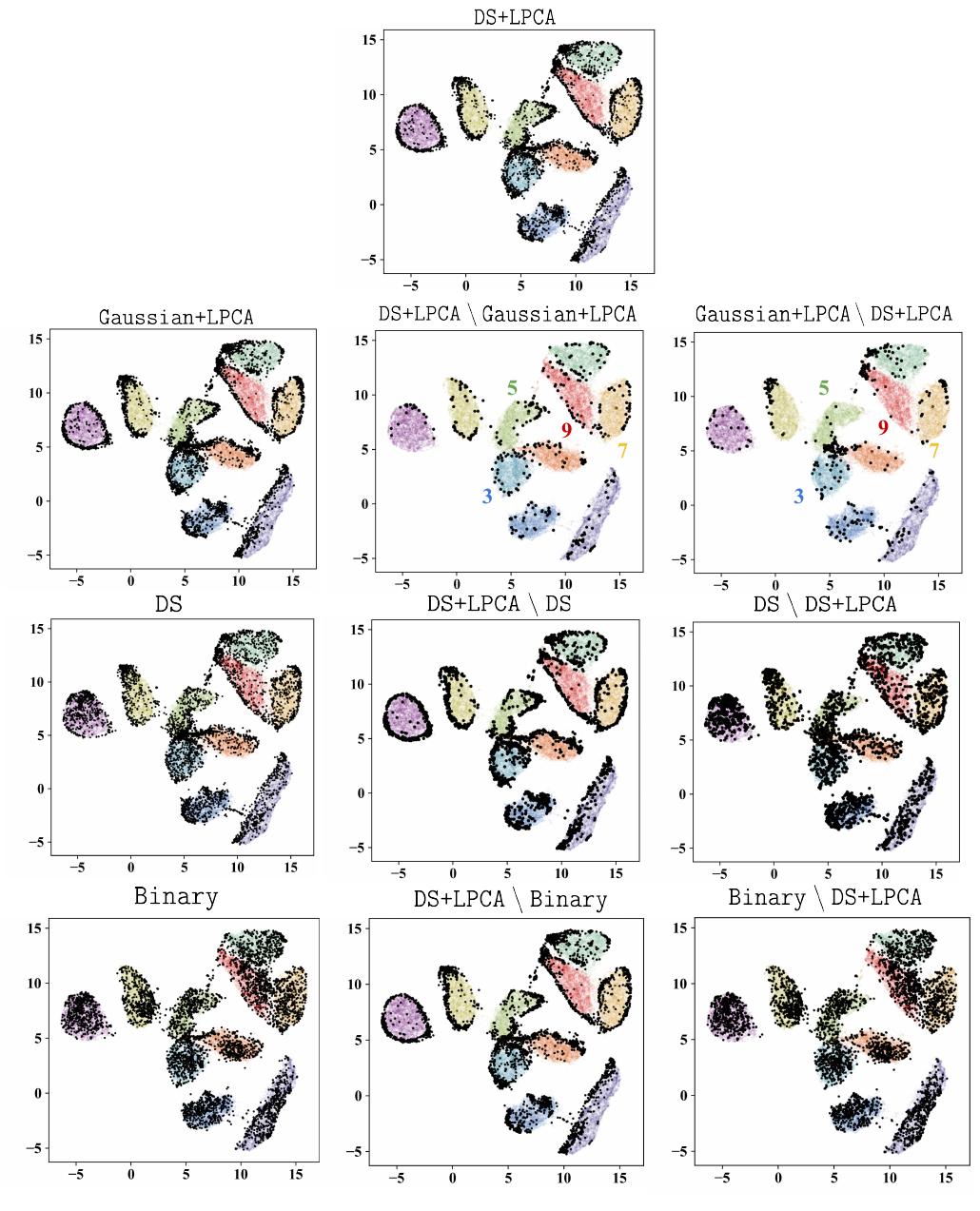}
    \end{tabular}
    \captionsetup{width=\linewidth}
    \vspace{0.2cm}
    \caption{The $10$-th percentile boundary $\widehat{\mathcal{B}}_{10}$ (as defined in Eq.~\cref{eq:Chatp2}-\cref{eq:Chatp}) is estimated from $5000$ raw $28 \times 28$ images of each MNIST digit using \textalgo{DS+LPCA}, \textalgo{Gaussian + LPCA}, \textalgo{DS} and \textalgo{Binary}. The visualizations are presented on a two-dimensional UMAP\cite{umap} embedding of the data, where points in $\widehat{\mathcal{B}}_{10}$ are highlighted in black. These plots show points on the estimated boundaries due to each method, the ones that lie on the estimated boundaries by two of the methods, and by one method but not the other.}
    \label{fig:mnist01}
\end{figure}

In this experiment, we randomly select $5000$ images of size $28 \times 28$ for each of the $10$ digits in the MNIST dataset, resulting in $10$ datasets, each containing $n=5000$ samples in $m=784$ dimensions. For each dataset, the median distance to the $128$th nearest neighbor of each image is used as the bandwidth $h$, and the local PCA is computed using $k_{\text{nn}} = 512$ and $d=2$. Using \textalgo{DS}, \textalgo{DS+LPCA}, \textalgo{Gaussian+LPCA} and \textalgo{Binary}, we compute the $10$th percentile boundary, $\widehat{\mathcal{B}}_{10}$, as defined in Eq.~\cref{eq:Chatp2}-\cref{eq:Chatp}. These results are visualized on a two-dimensional UMAP~\cite{umap} embedding of the entire dataset in~\Cref{fig:mnist01}. It is important to note that the boundary estimation is performed on the raw datasets, not the UMAP embeddings. Since we expect the boundary points to be localized and perhaps supported on curve segments, their visualization via UMAP can provide a qualitative assessment of the various boundary detection algorithms. It is also worth noting that certain points which seem to lie on the boundary of the clusters produced within the UMAP embedding in~\Cref{fig:mnist01} are not detected as boundary points by either \textalgo{DS+LPCA} or \textalgo{Gaussian+LPCA}.

Because there is no ground-truth boundary, we cannot evaluate accuracy directly as we did in previous sections. Instead, we are going use the Jaccard index to quantify the relative, qualitative agreement between the different methods. Treating \textalgo{DS+LPCA} as a reference, we find that the boundaries estimated by \textalgo{DS} and \textalgo{Binary} diverge significantly, yielding Jaccard indices of $0.257$ and $0.083$, respectively. It is visually evident from~\Cref{fig:mnist01} that boundary points detected by \textalgo{DS} and \textalgo{Binary} are highly scattered, inconsistent without our expectation of localized boundaries. In contrast, the boundaries computed using \textalgo{DS+LPCA} and \textalgo{Gaussian+LPCA} are concentrated along curve segments, and have a significant overlap with a Jaccard index of $0.847$. However, there are subtle differences. As shown in~\Cref{fig:mnist01}, points that lie on the boundary estimated by \textalgo{DS+LPCA} but not on the one by \textalgo{Gaussian+LPCA} tend to be more dispersed along the boundary, while \textalgo{Gaussian+LPCA} identifies points that are more concentrated in tightly packed regions resembling corners in the embedding, as is evident for digits $3$, $5$, $7$, and $9$. This suggests that boundary estimate by \textalgo{Gaussian+LPCA} is more skewed toward high-curvature, corner-like regions, in comparison to \textalgo{DS+LPCA}.

\begin{figure}
    \centering
    \includegraphics[width=0.975\linewidth]{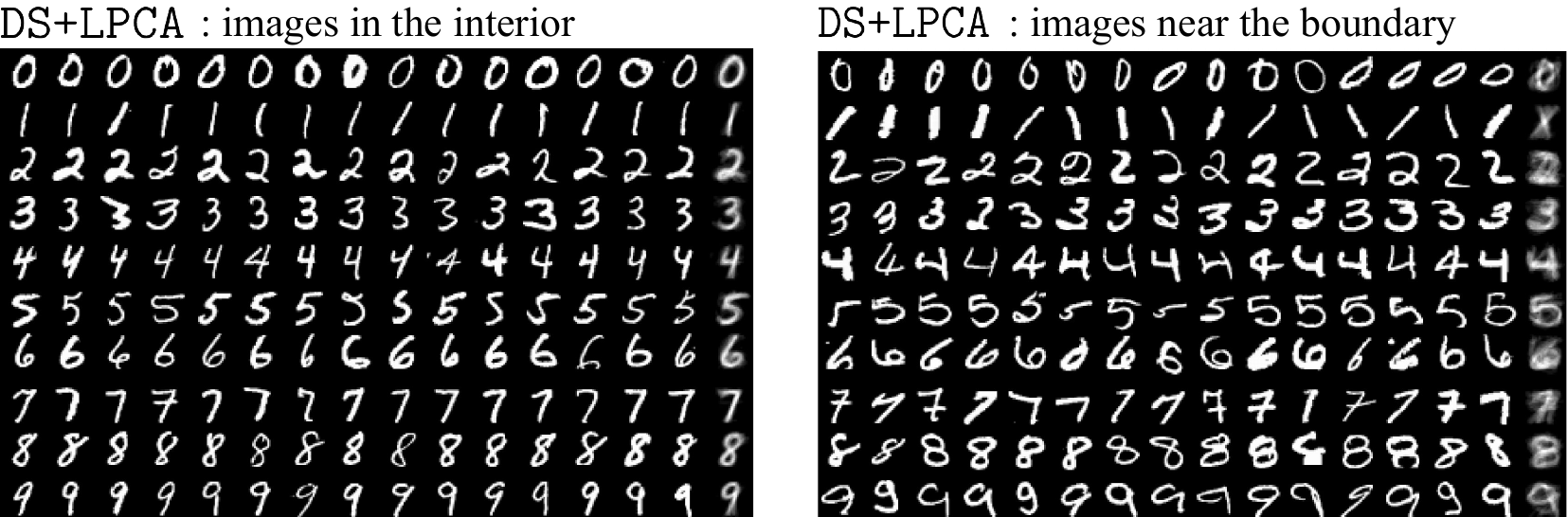}
    \vspace{0.2cm}
    \captionsetup{width=\linewidth}
    \caption{A random selection of 15 images per digit is shown, including those from the interior of the MNIST dataset (outside $\widehat{\mathcal{B}}_{95}$) and near the boundary (inside $\widehat{\mathcal{B}}_{5}$), as estimated using \textalgo{DS + LPCA}. The final column displays the average of the 15 images.}
    \label{fig:mnist2}
\end{figure}

Finally, to examine the differences between images near and far from the boundary, we randomly sampled $15$ images from $\widehat{\mathcal{B}}_{5}$ and outside $\widehat{\mathcal{B}}_{95}$, as estimated using our \textalgo{DS+LPCA}. The former set represents images near the boundary, while the latter represents images in the interior. These images are visualized in~\Cref{fig:mnist2}. Clearly, images within the interior exhibit less variation, whereas those near the boundary display greater variation. The contrast is evident in the average of the $15$ images per digit, displayed in the last columns of~\Cref{fig:mnist2}. This highlights that the boundary points capture more diverse and nuanced features of the dataset,  emphasizing the importance of accurately identifying them. Furthermore, this insight opens new avenues for exploring the computational benefits and generalization performance of machine learning models trained solely on boundary points rather than on the entire dataset. 

%% file: sections/conc.tex
In this work, we developed a robust framework for detecting boundary points in data sampled from noisy manifolds with boundary. To this end, we extended recent results on the doubly stochastic scaling of Gaussian kernels on closed manifolds~\cite{marshall2019manifold,landa2021doubly,landa2023robust,cheng2024bi} to the setting of manifolds with boundary. In particular, we derived a characterization of the scaling factors that transform a Gaussian kernel into a doubly stochastic kernel in the setting of a manifold with boundary. Building on this characterization, we proposed a BDE that combines the doubly stochastic kernel with local PCA, and utilized our characterization of the scaling factors to establish its convergence.

Numerical experiments demonstrate that the proposed BDE substantially outperforms standard Gaussian and binary kernel–based methods~\cite{berry2017density,calder2022boundary} in detecting boundary points, especially in high-noise regimes (\Cref{fig:jaccard_annulus,fig:jaccard_torus,fig:bde}). Moreover, we observe that the detected boundary points capture a broader range of salient features within the data (\Cref{fig:mnist2}).

We further introduced a two-parameter family of kernel density estimators that extends the KDE proposed in~\cite{landa2023robust} for closed manifolds to manifolds with boundary by incorporating a boundary correction term. Our experiments (\Cref{fig:kde_annulus,fig:kde_torus}) show that the proposed estimators yield improved density estimates compared to standard Gaussian KDEs with boundary correction.

Overall, this work highlights the importance of incorporating both doubly stochastic scaling and local PCA to achieve robust performance. As future work, we plan to investigate the trade-offs associated with training models using only boundary points versus the full dataset, with focus on on predictive performance, computational efficiency, and generalization.

%% file: sections/proofs.tex
The following definition and lemma are needed for the proofs that follow. Let $\Pi_{\mathbf{x}}:\mathcal{M} \mapsto T_{\mathbf{x}}\mathcal{M}$ be the projector onto the tangent space at $\mathbf{x}$. Define
\begin{equation}
    \mathcal{P}_{\mathbf{x}}(\mathbf{y}) \coloneqq  \Pi_{\mathbf{x}}(\mathbf{y}) - \mathbf{x}. \label{eq:Px_of_y}
\end{equation}
\begin{lemma}
\label{lemma:Px_eqs}
For $\mathbf{v}_{\mathbf{x}} \in T_{\mathbf{x}}\mathcal{M}$ and $\partial_{\mathbf{v}_{\mathbf{x}}} = \partial/\partial \mathbf{v}_{\mathbf{x}}$,
\begin{equation}
    \mathcal{P}_{\mathbf{x}}(\mathbf{x}) = 0,\ \partial_{\mathbf{v}_{\mathbf{x}}}\mathcal{P}_{\mathbf{x}}(\mathbf{x}) = \mathbf{v}_{\mathbf{x}}\  \text{ and }\ \partial_{\mathbf{v}_{\mathbf{x}}}^2\mathcal{P}_{\mathbf{x}}(\mathbf{x}) = 0.
\end{equation}
\end{lemma}

\begin{proof}
Let $(\partial_i)_1^d$ be an orthonormal basis of $T_{\mathbf{x}}\mathcal{M}$. Then there exist $a^{i} \in \mathbb{R}$ and $b^{i}:\mathcal{M} \mapsto \mathbb{R}$ for all $i \in [1,d]$ such that $\mathbf{v}_{\mathbf{x}} = a^{i}\partial_i$ (using Einstein notation) and $\mathcal{P}_{\mathbf{x}}(\mathbf{y}) = b^{i}(\mathbf{y})\partial_i$. Moreover, $b^{i}(\mathbf{x}) = 0$ for all $i \in [1,d]$,
\begin{align}
    \partial_{j}b^{i}(\mathbf{y}) &= \delta_{j}^{i} + \mathcal{O}(\left\|\mathbf{x}-\mathbf{y}\right\|^2) \text{and}\\
    \partial_{k}\partial_{j}b^{i}(\mathbf{y}) &= \mathcal{O}(\left\|\mathbf{x}-\mathbf{y}\right\|).
\end{align}
Thus, $\mathcal{P}_{\mathbf{x}}(\mathbf{x}) = 0$ is trivial. Then
\begin{align}
    \partial_{\mathbf{v}_{\mathbf{x}}}\mathcal{P}_{\mathbf{x}}(\mathbf{x}) &= a^{j}\partial_{j}b^{i}(\mathbf{x}) \partial_{i} = a^{j}\delta_{j}^{i}\partial_{i} = a^{i}\partial_{i} = \mathbf{v}_{\mathbf{x}},\\
    \partial_{\mathbf{v}_{\mathbf{x}}}^2\mathcal{P}_{\mathbf{x}}(\mathbf{x}) &= \mathbf{v}_{\mathbf{x}}(\mathbf{v}_{\mathbf{x}} \mathcal{P}_{\mathbf{x}})(\mathbf{x}) =  \mathbf{v}_{\mathbf{x}} (a^{j}\partial_{j}b^{i})(\mathbf{x}) \partial_{i} = a^{k}a^{j}\partial_{k}\partial_{j}b^{i}(\mathbf{x})\partial_{i} = 0
\end{align}
\end{proof}

\begin{proof}[Proof of~\Cref{thm:rho_hx_char}]
For brevity, we use the following shorthand notation for the proof: $\rho \equiv \rho_h(\mathbf{x})$, $q \equiv q(\mathbf{x})$, $H \equiv H(\mathbf{x})$, $\kappa_{\pm} \equiv \kappa_{\pm}(\mathbf{x})$, $A \equiv A_h(\mathbf{x})$, $B \equiv B_h(\mathbf{x})$, $m_i \equiv m_h^{(i)}(\mathbf{x})$ for $i=1,2$ and $3$, and $\widetilde{m}_2 = \frac{\pi^{d/2}}{2}$.

For part (I), take the derivative of Eq.~\cref{eq:rho_constraint} with respect to $\mathbf{x}$ in the direction of $\mathbf{v}_{\mathbf{x}} \in T_{\mathbf{x}}\mathcal{M}$,
\begin{equation}
    \frac{1}{h^d}\int_{\mathcal{M}}(k_{h}(\mathbf{x},\mathbf{y})\partial_{\mathbf{v}_{\mathbf{x}}}\rho_h(\mathbf{x}) + \rho_{h}(\mathbf{x})\ \partial_{\mathbf{v}_{\mathbf{x}}}k_{h}(\mathbf{x},\mathbf{y})) \rho_{h}(\mathbf{y}) q(\mathbf{y})dV(\mathbf{y}) = 0. \label{eq:partial_eta_rho_h_constraint}
\end{equation}

For the first term, using~\Cref{lem:integral_khf} and Eq.~\cref{eq:rho_constraint}, we obtain:
\begin{align}
    &\frac{1}{h^d}\int_{\mathcal{M}}k_{h}(\mathbf{x},\mathbf{y})\partial_{\mathbf{v}_{\mathbf{x}}}\rho_h(\mathbf{x}) \rho_{h}(\mathbf{y}) q(\mathbf{y})dV(\mathbf{y})\\
    &\hspace{1cm}= \pi^{d/2}\frac{\partial_{\mathbf{v}_{\mathbf{x}}}\rho_h(\mathbf{x})}{\rho_h(\mathbf{x})} \frac{1}{\pi^{d/2}h^d}\int_{\mathcal{M}}k_{h}(\mathbf{x},\mathbf{y})\rho_h(\mathbf{x}) \rho_{h}(\mathbf{y}) q(\mathbf{y})dV(\mathbf{y})\\
    &\hspace{1cm}= \pi^{d/2}\frac{\partial_{\mathbf{v}_{\mathbf{x}}}\rho_h(\mathbf{x})}{\rho_h(\mathbf{x})}. \label{eq:dvxrho1}
\end{align}
Then, for the second term, using $\partial_{\mathbf{v}_{\mathbf{x}}}k_{h}(\mathbf{x},\mathbf{y}) = 2k_{h}(\mathbf{x},\mathbf{y})\mathbf{v}_{\mathbf{x}}^T\mathcal{P}_{\mathbf{x}}(\mathbf{y})/h^2$ followed by~\Cref{lem:integral_khf,lemma:Px_eqs} and the following facts
\begin{align}
    \Delta(f \mathbf{v}_{\mathbf{x}}^T\mathcal{P}_{\mathbf{x}})(\mathbf{x}) &= \mathbf{v}_{\mathbf{x}}^T\mathcal{P}_{\mathbf{x}}(\mathbf{x})\Delta f(\mathbf{x}) + f(\mathbf{x})\Delta(\mathbf{v}_{\mathbf{x}}^T\mathcal{P}_{\mathbf{x}})(\mathbf{x}) + 2 \partial_{\mathbf{v}_{\mathbf{x}}}f(\mathbf{x})\\
    &= 0 + 0 + 2 \partial_{\mathbf{v}_{\mathbf{x}}}f(\mathbf{x})\\
    \partial_{\boldsymbol{\eta}_{\mathbf{x}}}^2(f \mathbf{v}_{\mathbf{x}}^T\mathcal{P}_{\mathbf{x}})(\mathbf{x}) &= 0 + 0 + 2\partial_{\boldsymbol{\eta}_{\mathbf{x}}}f(\mathbf{x})\mathbf{v}_{\mathbf{x}}^T\boldsymbol{\eta}_{\mathbf{x}}\\
    \partial_{\boldsymbol{\eta}_{\mathbf{x}}}(f \mathbf{v}_{\mathbf{x}}^T\mathcal{P}_{\mathbf{x}})(\mathbf{x}) &= f(\mathbf{x}) \mathbf{v}_{\mathbf{x}}^T\boldsymbol{\eta}_{\mathbf{x}},
\end{align}
we obtain the following equations:
\begin{align}
    &\frac{2\rho_{h}(\mathbf{x})}{h^{d+2}}\int_{\mathcal{M}} k_{h}(\mathbf{x},\mathbf{y})\mathbf{v}_{\mathbf{x}}^T\mathcal{P}_{\mathbf{x}}(\mathbf{y})  \rho_{h}(\mathbf{y}) q(\mathbf{y})dV(\mathbf{y})\\
    &\hspace{0.5cm}=\frac{2 m_h^{(1)}(\mathbf{x})\rho_{h}(\mathbf{x})^2q(\mathbf{x})\mathbf{v}_{\mathbf{x}}^T\boldsymbol{\eta}_{\mathbf{x}}}{h} + \rho_h(\mathbf{x}) \left\{2\widetilde{m}_2\partial_{\mathbf{v}_{\mathbf{x}}}(\rho_h q)(\mathbf{x}) +\right.\\
    &\hspace{1.25cm}  \left. 2(m_h^{(2)}(\mathbf{x})-\widetilde{m}_2)\partial_{\boldsymbol{\eta}_{\mathbf{x}}}(\rho_h q)(\mathbf{x})\mathbf{v}_{\mathbf{x}}^T\boldsymbol{\eta}_{\mathbf{x}} - (d-1)m_h^{(2)}(\mathbf{x})H(\mathbf{x})\rho_h(\mathbf{x})q(\mathbf{x}) \mathbf{v}_{\mathbf{x}}^T\boldsymbol{\eta}_{\mathbf{x}}\right\} + \mathcal{O}(h) \label{eq:dvxrho2}
\end{align}
By substituting $\mathbf{v}_{\mathbf{x}} = \boldsymbol{\eta}_{\mathbf{x}}$ in Eq.~\cref{eq:dvxrho1} and Eq.~\cref{eq:dvxrho2}, and equating the sum to zero, we obtain
\begin{align}
    \pi^{d/2}\frac{\partial_{\boldsymbol{\eta}_{\mathbf{x}}}\rho}{\rho} + \frac{2 m_1 \rho^2q }{h} + \rho\left\{2m_2 \partial_{\boldsymbol{\eta}_{\mathbf{x}}}(\rho q) - (d-1)m_2H\rho q\right\} + \mathcal{O}(h) &= 0\\
    (\pi^{d/2}+2m_2\rho^2 q)h\partial_{\boldsymbol{\eta}_{\mathbf{x}}}\rho + 2\rho^3(m_1q + hm_2\partial_{\boldsymbol{\eta}_{\mathbf{x}}}q) -(d-1)hm_2H\rho^3 q + \mathcal{O}(h^2) &= 0\\
    (\pi^{d/2}+2m_2\rho^2 q)h\partial_{\boldsymbol{\eta}_{\mathbf{x}}}\rho + 2\rho^3q(m_1 + hm_2 \kappa_{-}) + \mathcal{O}(h^2) &= 0\label{eq:partial_eta_x_rho_short}
\end{align}
Similarly, substituting $\mathbf{v}_{\mathbf{x}}^T\boldsymbol{\eta}_{\mathbf{x}} = 0$ in Eq.~\cref{eq:dvxrho1} and Eq.~\cref{eq:dvxrho2}, and equating the sum to zero, we obtain
\begin{equation}
    \pi^{d/2}\frac{\partial_{\mathbf{v}_{\mathbf{x}}}\rho}{\rho} + 2\widetilde{m}_2 \rho\left(q\partial_{\mathbf{v}_{\mathbf{x}}}\rho + \rho\partial_{\mathbf{v}_{\mathbf{x}}}q\right) + \mathcal{O}(h) = 0.
\end{equation}
Simplifying the above and substituting the value of $\widetilde{m}_2 = \pi^{d/2}/2$, we obtain
\begin{equation}
    (1+\rho^2 q)\partial_{\mathbf{v}_{\mathbf{x}}}\rho + \rho^3\partial_{\mathbf{v}_{\mathbf{x}}}q + \mathcal{O}(h) = 0.
\end{equation}

For part (II), we expand Eq.~\cref{eq:rho_constraint} using \Cref{lem:integral_khf},
\begin{align}
    \pi^{d/2} &= \rho\left\{m_0 \rho q + hm_1\left(\partial_{\boldsymbol{\eta}_{\mathbf{x}}}(\rho q)+ \frac{d-1}{2}H\rho q\right)\right\} + \mathcal{O}(h^2)\\
    &= \rho\left\{\left(m_0 + hm_1\frac{d-1}{2}H\right)\rho q + hm_1\partial_{\boldsymbol{\eta}_{\mathbf{x}}}(\rho_hq)\right\} + \mathcal{O}(h^2)\\
    &= \rho q\left\{(m_0 + hm_1 \kappa_{+})\rho + hm_1\partial_{\boldsymbol{\eta}_{\mathbf{x}}}\rho\right\} + \mathcal{O}(h^2) \label{eq:eq1}
\end{align}
Rearranging the above equation, we obtain 
\begin{equation}
    hm_1\rho q\partial_{\boldsymbol{\eta}_{\mathbf{x}}}\rho = \pi^{d/2} - (m_0 + hm_1\kappa_{+})\rho^2q + \mathcal{O}(h^2).\label{eq:partial_eta_rho_h_short}
\end{equation}
Multiplying Eq.~\cref{eq:partial_eta_x_rho_short} by $m_1\rho q$,
\begin{equation}
    (\pi^{d/2}+2m_2\rho^2 q)hm_1\rho q\partial_{\boldsymbol{\eta}_{\mathbf{x}}}\rho + 2m_1\rho^4q^2(m_1 + hm_2 \kappa_{-}) + \mathcal{O}(h^2) = 0
\end{equation}
and then substituting Eq.~\cref{eq:partial_eta_rho_h_short} results in
\begin{equation}
    A\rho^4 q^2- B\rho^2 q  + \pi^d = \mathcal{O}(h^2)
\end{equation}
where, using the definition of $\kappa_{+}$ in the statement of the theorem,
\begin{align}
    A &= 2\left(m_1^2 + hm_1m_2\kappa_{-} - m_2(m_0+hm_1\kappa_{+})\right)\\
    &= 2\left(m_1^2 - m_2\left(m_0 + h(d-1)m_1H\right)\right)
\end{align}
and, using the definition of $m_2$ in Eq.~\cref{eq:m2},
\begin{equation}
    B = \pi^{d/2}(m_0 + hm_1\kappa_{+} - 2m_2) = \pi^{d/2}m_1\left(h\kappa_{+} - \frac{2b_{\mathbf{x}}}{h}\right).
\end{equation}
\end{proof}

\begin{proof}[Proof of~\Cref{corollary:closed_manifold}]
The proof follows from Eq.~\cref{eq:partial_eta_rho_h_short} by using the fact $b_{\mathbf{x}} = \infty$ for every point on a closed manifold, and as a result of substituting it in Eq.~\cref{eq:m0,eq:m1}, we obtain $m^{(0)}_{h}(\mathbf{x}) = \pi^{d/2}$ and $m^{(1)}_{h}(\mathbf{x}) = 0$.
\end{proof}

\begin{proof}[Proof of~\Cref{corollary:rho_h_H_zero}]
We adopt the same notation as in the proof of Theorem~\ref{thm:rho_hx_char}.
Due to the assumption $H = 0$, we have
\begin{equation}
    A = 2(m_1^2 - m_2m_0) \in \left[-\frac{\pi^{d-1}}{4}\left(\pi-1\right), -\pi^{d}\right).
\end{equation}
Consequently, $B^2 - 4\pi^d A > 0$ and therefore its square root is well defined. For convenience, define $C_{\pm} \coloneqq \frac{B \pm \sqrt{B^2 - 4\pi^dA}}{2A}$ and note that $C_h(\mathbf{x}) = C_{-}$. Then, Eq.~\cref{eq:rho_h_quadratic} reduces to
\begin{equation}
    (\rho^2q - C_{-})(\rho^2q + C_{+}) = \mathcal{O}(h^2).
\end{equation}
Since $\rho > 0$, $C_{+} < 0$ and
\begin{equation}
    C_{+} - C_{-} = \frac{2\sqrt{B^2 - 4\pi^dA}}{2A} \geq \frac{2\pi^{d/2}}{\sqrt{|A|}} \geq 2,
\end{equation}
therefore $\rho^2q = C_{-} + \mathcal{O}(h^2)$. Since $\mathcal{M}$ is compact therefore $C_{-}$ is bounded and thus $\rho^2q = C_{-}(1 + \mathcal{O}(h^2))$. Finally, if $\partial_{\boldsymbol{\eta}_{\mathbf{x}}}q$ is also zero then $\kappa_+ = 0$ by definition and the result follows.
\end{proof}

\begin{proof}[Proof of~\Cref{corollary:rho_h_boundary}]
Setting $b_{\mathbf{x}} = 0$, we obtain $m^{(0)}_{h}(\mathbf{x}) = \frac{\pi^{d/2}}{2}$, $m^{(1)}_{h}(\mathbf{x}) = -\frac{\pi^{(d-1)/2}}{2}$ and $m^{(2)}_{h}(\mathbf{x}) = \frac{\pi^{d/2}}{4}$.
\end{proof}

\begin{proof}[Proof of~\Cref{corollary:rho_h_first_order_char}]
This follows trivially from Theorem~\ref{thm:rho_hx_char} by taking away the terms involving $hH(\mathbf{x})$ and $h\kappa_{+}(\mathbf{x})$ and solving the quadratic in $\rho_h^2q$.
\end{proof}

\begin{proof}[Proof of~\Cref{theorem:Enu}]
Using assumption \refassump{assump:tangent}{(A6)} followed by \refassump{assump:ortho_noise}{(A7)}, we have
\begin{equation}
    \boldsymbol{\nu}_i = \frac{1}{n-1}\sum_{j=1}^{n}\mathbf{W}_{ij}\boldsymbol{\mathcal{U}}_{i}^T(\mathbf{x}_j + \boldsymbol{\varepsilon}(\mathbf{x}_j) - (\mathbf{x}_i + \boldsymbol{\varepsilon}(\mathbf{x}_i))) = \frac{1}{n-1}\sum_{j=1}^{n}\mathbf{W}_{ij}\boldsymbol{\mathcal{U}}_{i}^T(\mathbf{x}_j - \mathbf{x}_i).
\end{equation}
Using Eq.~\cref{eq:Wij_result} and the fact that $\left\|\mathbf{x}_k\right\| \leq 1$ and $\boldsymbol{\mathcal{U}}_{i}^T\boldsymbol{\mathcal{U}}_{i} = \mathbb{I}_d$,
\begin{equation}
    \boldsymbol{\nu}_i = \frac{\rho_h(\mathbf{x}_i)}{\pi^{d/2}h^d}\ \boldsymbol{\mathcal{U}}_{i}^T\left(\int k_h(\mathbf{x}_i, \mathbf{y})\rho_h(\mathbf{y}) (\mathbf{y}-\mathbf{x}_i)q(\mathbf{y})dV(\mathbf{y})\right) + \mathcal{O}^{(h)}_{m,n}(g(m,n)).
\end{equation}
From here on, for brevity, we use the shorthand notation as in the proof of Theorem~\ref{thm:rho_hx_char} and use $\mathbf{x}$ in place of $\mathbf{x}_i$.
Due to \Cref{lem:integral_khf} and the product rule for second order derivatives,
{\small
\begin{equation}
    \frac{1}{h^d}\int k_h(\mathbf{x}, \mathbf{y}) \rho_h(\mathbf{y})(\mathbf{y}-\mathbf{x}) q(\mathbf{y})dV(\mathbf{y}) = hm_1\rho q \boldsymbol{\eta}_{\mathbf{x}} + h^2\left(\widetilde{m}_2\nabla (\rho q) + \left(m_2 - \widetilde{m}_2\right)\partial_{\boldsymbol{\eta}_{\mathbf{x}}}(\rho q)\boldsymbol{\eta}_{\mathbf{x}}\right) + \mathcal{O}\left(h^2 \right)\boldsymbol{\eta}_{\mathbf{x}}
\end{equation}}%
and from~\Cref{thm:rho_hx_char} we have,
\begin{equation}
    h\partial_{\boldsymbol{\eta}_{\mathbf{x}}}\rho = -\frac{2m_1\rho^3 q}{\pi^{d/2}+2m_2\rho^2q} + \mathcal{O}(h)
\end{equation}
and, if $\mathbf{v}_\mathbf{x} \perp \boldsymbol{\eta}_{\mathbf{x}}$ then
\begin{equation}
    \partial_{\mathbf{v}_\mathbf{x}}(\rho q) = (\rho\partial_{\mathbf{v}_\mathbf{x}}q + q\partial_{\mathbf{v}_\mathbf{x}}\rho) = \left(\rho\partial_{\mathbf{v}_\mathbf{x}}q - \frac{\rho^3q\partial_{\mathbf{v}_\mathbf{x}}q}{1+\rho^2q}\right) + \mathcal{O}(h)= \frac{\rho\partial_{\mathbf{v}_\mathbf{x}}q}{1+\rho^2q} + \mathcal{O}(h),
\end{equation}
where the constants in $\mathcal{O}\left(h^2\right)$ and $\mathcal{O}\left(h\right)$, respectively, can be bounded uniformly over $\mathbf{x} \in \mathcal{M}$ (this holds throughout the proof).
Let $\boldsymbol{\eta}_{\mathbf{x}_{\perp}}$ represent the subspace orthogonal to $\boldsymbol{\eta}_{\mathbf{x}}$, then $\partial_{\boldsymbol{\eta}_{\mathbf{x}_{\perp}}}(\rho q)\boldsymbol{\eta}_{\mathbf{x}_{\perp}} = \nabla (\rho q) - \partial_{\boldsymbol{\eta}_{\mathbf{x}}}(\rho q)\boldsymbol{\eta}_{\mathbf{x}}$. Using the above equations, we obtain
\begin{align}
    &hm_1\rho q \boldsymbol{\eta}_{\mathbf{x}} + h^2(\widetilde{m}_2\partial_{\boldsymbol{\eta}_{\mathbf{x}_{\perp}}}(\rho q)\boldsymbol{\eta}_{\mathbf{x}_{\perp}} + m_2 \partial_{\boldsymbol{\eta}_{\mathbf{x}}}(\rho q)\boldsymbol{\eta}_{\mathbf{x}}) + \mathcal{O}\left(h^2\right) \boldsymbol{\eta}_{\mathbf{x}}\\
    &= hm_1\rho q \boldsymbol{\eta}_{\mathbf{x}} + h^2(\widetilde{m}_2\partial_{\boldsymbol{\eta}_{\mathbf{x}_{\perp}}}(\rho q)\boldsymbol{\eta}_{\mathbf{x}_{\perp}} + m_2 (\rho \partial_{\boldsymbol{\eta}_{\mathbf{x}}}q + q\partial_{\boldsymbol{\eta}_{\mathbf{x}}}\rho)\boldsymbol{\eta}_{\mathbf{x}}) + \mathcal{O}\left(h^2\right)\boldsymbol{\eta}_{\mathbf{x}}\\
    &= h\left(m_1\rho q  - \frac{2m_1m_2\rho^3q^2}{\pi^{d/2}+2m_2 \rho^2q}\right)\boldsymbol{\eta}_{\mathbf{x}}+ h^2 \widetilde{m}_2\partial_{\boldsymbol{\eta}_{\mathbf{x}_{\perp}}}(\rho q)\boldsymbol{\eta}_{\mathbf{x}_{\perp}}  + \mathcal{O}\left(h^2\right)\boldsymbol{\eta}_{\mathbf{x}}\\
    &= h\left(\frac{\pi^{d/2}m_1\rho q}{\pi^{d/2}+2m_2 \rho^2q}\right)\boldsymbol{\eta}_{\mathbf{x}}+ h^2\widetilde{m}_2\partial_{\boldsymbol{\eta}_{\mathbf{x}_{\perp}}}(\rho q)\boldsymbol{\eta}_{\mathbf{x}_{\perp}}  + \mathcal{O}\left(h^2\right)\boldsymbol{\eta}_{\mathbf{x}}\\
    &h\left(\frac{\pi^{d/2}m_1\rho q}{\pi^{d/2}+2m_2 \rho^2q}\right)\boldsymbol{\eta}_{\mathbf{x}}+ \frac{h^2\pi^{d/2}\rho}{2}\frac{\partial_{\boldsymbol{\eta}_{\mathbf{x}_{\perp}}}q}{1+\rho^2q}\boldsymbol{\eta}_{\mathbf{x}_{\perp}}   + \mathcal{O}\left(h^2\right)\boldsymbol{\eta}_{\mathbf{x}}\\
    &= \frac{h}{\rho}\left(\frac{\pi^{d/2}m_1\zeta}{\pi^{d/2}+2m_2 \zeta}\right)\boldsymbol{\eta}_{\mathbf{x}}  + \mathcal{O}\left(h^2\right)\boldsymbol{\eta}_{\mathbf{x}} + \mathcal{O}\left(h^2\right)\boldsymbol{\eta}_{\mathbf{x}_{\perp}}.
\end{align}
Using the definition of $\beta_h$ in the statement of the theorem and the fact that $\boldsymbol{\mathcal{U}}_{i}^T\boldsymbol{\eta}_{\mathbf{x}_i} = \boldsymbol{\eta}_{\mathbf{x}_i}$, we obtain the result.
\end{proof}

\begin{proof}[Proof of~\Cref{corollary:Enu}]
The proof closely follows that of the previous theorem. Below, we detail the specific equations modified by the newly introduced assumptions. As before, we adopt the shorthand notation and use $\mathbf{x}$ in place of $\mathbf{x}_i$. 

Given $H(\mathbf{x}) = 0$, applying \Cref{lem:integral_khf} and the product rule for second order derivatives yields,
{\small
\begin{equation}
    \frac{1}{h^d}\int k_h(\mathbf{x}, \mathbf{y}) \rho_h(\mathbf{y})(\mathbf{y}-\mathbf{x}) q(\mathbf{y})dV(\mathbf{y}) = hm_1\rho q \boldsymbol{\eta}_{\mathbf{x}} + h^2\left(\widetilde{m}_2\nabla (\rho q) + \left(m_2 - \widetilde{m}_2\right)\partial_{\boldsymbol{\eta}_{\mathbf{x}}}(\rho q)\boldsymbol{\eta}_{\mathbf{x}}\right) + \mathcal{O}\left(h^3 \right)
\end{equation}}%
Moreover, because $q$ is uniform, we have $\kappa_{-}(\mathbf{x}) = 0$. Consequently~\Cref{thm:rho_hx_char} implies:
\begin{equation}
    h\partial_{\boldsymbol{\eta}_{\mathbf{x}}}\rho = -\frac{2m_1\rho^3 q}{\pi^{d/2}+2m_2\rho^2q} + \mathcal{O}(h^2).
\end{equation}
Furthermore, for any vector $\mathbf{v}_\mathbf{x} \perp \boldsymbol{\eta}_{\mathbf{x}}$, it follows that:
\begin{equation}
    \partial_{\mathbf{v}_\mathbf{x}}(\rho q) = (\rho\partial_{\mathbf{v}_\mathbf{x}}q + q\partial_{\mathbf{v}_\mathbf{x}}\rho) = \left(\rho\partial_{\mathbf{v}_\mathbf{x}}q - \frac{\rho^3q\partial_{\mathbf{v}_\mathbf{x}}q}{1+\rho^2q}\right) + \mathcal{O}(h)=  \mathcal{O}(h).
\end{equation}
The result follows by substituting the updated equations in the previous proof.
\end{proof}

\begin{proof}[Proof of~\Cref{theorem:Eq}]
Using Eq.~(\ref{eq:Wij_result}, \ref{eq:g^h(X)}, \ref{eq:maxXi}),
\begin{align}
    \widehat{q}_i^{1-s} &= \frac{((n-1)h^d)^{s-1}\pi^{ds/2}s^{d/2}}{(n-1)^s\pi^{ds/2}h^{ds}} \frac{1}{(m^{(0)}_{h/\sqrt{s}}(\mathbf{x}_i)\zeta_h(\mathbf{x}_i)^{s})^{u}} \times \\
    &\hspace{2cm} \sum_{j=1}^{n} \frac{\rho_h(\mathbf{x}_i)^s k_h(\mathbf{x}_i, \mathbf{x}_j)^s\rho_h(\mathbf{x}_j)^s}{(m^{(0)}_{h/\sqrt{s}}(\mathbf{x}_j)\zeta_h(\mathbf{x}_j)^{s})^{1-u}}(1 + \mathcal{O}^{(h)}_{m,n}(g(m,n)))^s\\
    &= \frac{1}{(n-1)(h/\sqrt{s})^d}  \frac{1}{(m^{(0)}_{h/\sqrt{s}}(\mathbf{x}_i)\zeta_h(\mathbf{x}_i)^{s})^{u}} \times\\
    &\hspace{2cm} \sum_{j=1}^{n} \frac{\rho_h(\mathbf{x}_i)^s k_h(\mathbf{x}_i, \mathbf{x}_j)^s\rho_h(\mathbf{x}_j)^s}{(m^{(0)}_{h/\sqrt{s}}(\mathbf{x}_j)\zeta_h(\mathbf{x}_j)^{s})^{1-u}}(1 + \mathcal{O}^{(h)}_{m,n}(g(m,n)))\\
    &= \left(\frac{1}{(h/\sqrt{s})^d}  \frac{\rho_h(\mathbf{x}_i)^s}{(m^{(0)}_{h/\sqrt{s}}(\mathbf{x}_i)\zeta_h(\mathbf{x}_i)^{s})^u}  \times \right.\\
    &\hspace{2cm} \left. \int  \frac{k_{h/\sqrt{s}}(\mathbf{x}_i, \mathbf{y})\rho_h(\mathbf{y})^sq(\mathbf{y})}{(m^{(0)}_{h/\sqrt{s}}(\mathbf{y})\zeta_h(\mathbf{y})^{s})^{1-u}}dV(\mathbf{y}) + \mathcal{O}^{(h)}_{m,n}\left(\sqrt{\frac{\log n}{n}}\right) \right)(1 + \mathcal{O}^{(h)}_{m,n}(g(m,n)))\\
    &= \left(\frac{1}{(h/\sqrt{s})^d}  \frac{\rho_h(\mathbf{x}_i)^s}{(m^{(0)}_{h/\sqrt{s}}(\mathbf{x}_i)\zeta_h(\mathbf{x}_i)^{s})^u}  \times \right.\\
    &\hspace{2cm} \left.\int  \frac{k_{h/\sqrt{s}}(\mathbf{x}_i, \mathbf{y})\rho_h(\mathbf{y})^sq(\mathbf{y})}{(m^{(0)}_{h/\sqrt{s}}(\mathbf{y})\zeta_h(\mathbf{y})^{s})^{1-u}}dV(\mathbf{y}) \right)(1 + \mathcal{O}^{(h)}_{m,n}(g(m,n))). \label{eq:temp1}
\end{align}
where we used the Hoeffding's inequality conditioned on $\mathbf{x}_i$ and Eq.~(\ref{eq:g^h(X)}) in the last two equations, combined with the facts that $\zeta_h$ is bounded away from zero by definition, and $\rho_h$ and $k_h$ are bounded due to compactness of $\mathcal{M}$. We then proceed by invoking Lemma~\ref{lem:integral_khf} followed by Corollary~\ref{corollary:rho_h_first_order_char},
\begin{align}
    &\frac{\rho_h(\mathbf{x}_i)^s}{(m^{(0)}_{h/\sqrt{s}}(\mathbf{x}_i)\zeta_h(\mathbf{x}_i)^{s})^u}\frac{1}{(h/\sqrt{s})^d} \int  \frac{k_{h/\sqrt{s}}(\mathbf{x}_i, \mathbf{y})\rho_h(\mathbf{y})^sq(\mathbf{y})}{(m^{(0)}_{h/\sqrt{s}}(\mathbf{y})\zeta_h(\mathbf{y})^{s})^{1-u}}dV(\mathbf{y})\\
    &\qquad=  \frac{\rho_h(\mathbf{x}_i)^s}{(m^{(0)}_{h/\sqrt{s}}(\mathbf{x}_i)\zeta_h(\mathbf{x}_i)^{s})^u}\frac{m^{(0)}_{h/\sqrt{s}}(\mathbf{x}_i)\rho_h(\mathbf{x}_i)^sq(\mathbf{x}_i)}{(m^{(0)}_{h/\sqrt{s}}(\mathbf{x}_i)\zeta_h(\mathbf{x}_i)^{s})^{1-u}} + \mathcal{O}\left(\frac{h}{\sqrt{s}}\right)\\
    &\qquad= \frac{\rho_h(\mathbf{x}_i)^{2s}}{\zeta_h(\mathbf{x}_i)^{s}}q(\mathbf{x}_i)+  \mathcal{O}\left(\frac{h}{\sqrt{s}}\right)\\
    &\qquad= q(\mathbf{x}_i)^{1-s} + \max\left\{s,\frac{1}{\sqrt{s}}\right\}\mathcal{O}(h).
\end{align}
The result follows by combining Eq.~(\ref{eq:temp1}) with the above equation, taking a power of $1/(1-s)$ and applying Eq.~(\ref{eq:g^h(X)}), while noting that $m^{(0)}_{h/\sqrt{s}}$ is bounded away from zero by definition, $q$ is bounded away from zero by assumption \refassump{assump:q}{(A2')} and $\mathcal{M}$ is compact.
\end{proof}